\documentclass[10pt,twocolumns]{IEEEtran}

\usepackage{cite}
\usepackage{amsmath,amssymb,amsfonts}
\usepackage{algorithmic}
\usepackage{textcomp}
\usepackage{color}
\usepackage{bm}
\usepackage{graphicx}
\usepackage{enumerate}
\usepackage{multicol}

\newtheorem{lemma}{Lemma}[section]
\newtheorem{assumption}{Assumption}[section]

\newtheorem{remark}{Remark}[section]

\newtheorem{theorem}{Theorem}[section]
\newtheorem{example}{Example}[section]

\numberwithin{equation}{section}

\title{
\huge{
Passivity-based Gradient-Play Dynamics for Distributed Generalized Nash Equilibrium Seeking
}
}

\author{Weijian~Li and Lacra Pavel
\thanks{W.~Li and L.~Pavel are with the Department of Electrical and Computer Engineering, University of Toronto, Toronto, ON, M5S 3G4, Canada. E-mails: weijian.li@utoronto.ca, pavel@control.utoronto.ca.
}
}

\begin{document}

\maketitle

\begin{abstract}
We consider seeking generalized Nash equilibria (GNE) for noncooperative games with  coupled nonlinear constraints over networks.
We first revisit a well-known gradient-play dynamics from a passivity-based perspective, and address that the strict monotonicity on pseudo-gradients is a critical assumption to ensure the exact convergence of the dynamics.
Then we propose two novel passivity-based gradient-play dynamics by introducing parallel feedforward compensators (PFCs) and output feedback compensators (OFCs).
We show that the proposed dynamics can reach exact GNEs  in merely monotone regimes if the PFCs are strictly passive or the OFCs are output strictly passive. 
Following that, resorting to passivity, we develop a unifying framework to generalize the gradient-play dynamics, and moreover, design a class of explicit passive-based dynamics with convergence guarantees.
In addition, we explore the relation between the proposed dynamics and some existing methods, and extend our results to partial-decision information settings.
\end{abstract}


\section{INTRODUCTION}
\label{sec:introduction}

In recent years, noncooperative games  with
coupled constraints have received much attention since they could model many decision-making problems in multi-agent systems \cite{hu2022distributed}.
Examples can be found in various areas, including 
demand-side management in smart grids \cite{saad2012game}, electric vehicle charging \cite{paccagnan2018nash},
wireless communication networks \cite{yin2011nash} and
social networks \cite{ghaderi2014opinion}.
The basic idea is that in a multi-agent system,
each player aims to minimize a local cost function depending on its own action as well as on the strategies of its opponents, and meanwhile, the shared constraints should be satisfied. A generalized Nash equilibrium (GNE) is a reasonable solution to such a problem, whereby no player can decrease its local cost by unilaterally changing its own decision \cite{facchinei2003finite, facchinei2010generalized}.
Hence, a flurry of research has been devoted to 
distributed GNE seeking, and various algorithms have been explored.
To name just a few, there are 
(projected) gradient-play dynamics, operator splitting approaches and payoff-based learning \cite{hu2022distributed, yi2019operator, tatarenko2018learning}.

One of the most studied methods is the gradient-play dynamics because it is easily implemented under full and partial-decision information settings.
Convergence of the dynamics to Nash equilibriums (NEs) was proved in \cite{flaam2002equilibrium, shamma2005dynamic}.
Then the dynamics was extended to NE seeking over networks  in \cite{gadjov2018passivity}, where each player only knew its neighbors' decision information.
Following that, a geometrically convergent discrete-time gradient-play algorithm was designed in \cite{tatarenko2020geometric}, 
and achieved faster rates by incorporating a Nesterov's accelerated scheme.
Resorting to a primal-dual framework, a fully distributed gradient-play dynamics was explored for noncooperative games with separable nonlinear coupled constraints in \cite{bianchi2021continuous}.
Similar ideas were used to seek GNEs for multi-cluster games with nonsmooth payoff functions, coupled  and set constraints in \cite{zeng2019generalized}.
However, all algorithms mentioned above require strict or strong monotonicity of the pseudo-gradients to ensure their convergence.
Unfortunately, the assumption fails in practical applications, including zero-sum games, Cournot games and resource allocation problems \cite{parise2014mean, scutari2014real, kannan2012distributed}.

Passivity, a powerful tool for the analysis and design of 
control systems, has had a profound impact on control theory. It has been employed to study the stability of (interconnected) dynamical systems, as well as to design controllers with stability guarantees \cite{sepulchre2012constructive, khalil2002nonlinear}.
Recently, the concept has been applied to (G)NE computation for noncooperative games \cite{pavel2022dissipativity, ye2023distributed, gao2020passivity, gadjov2022exact, gao2023second}.
On the one hand, it provides a better understanding of some existing algorithms by explaining why they can work under certain game settings. On the other hand, it guides us to develop some novel dynamics.
For instance, stable games and evolutionary dynamics were analyzed  in \cite{fox2013population} after they were modeled as passive dynamical systems.
In  \cite{gadjov2018passivity}, passivity was employed to prove the convergence of a gradient-play dynamics, while a passivity-based methodology was proposed in \cite{gao2020passivity} for the analysis and design of reinforcement learning dynamics in multi-agent finite games. 
In \cite{romano2020gne}, autonomous agents with passive dynamics were taken into consideration, and then, an inexact penalty method, reaching 
an $\epsilon$-GNE, was designed.
To relax the strict monotonicity assumption on pseudo-gradients, a heavy-anchor (HA) dynamics was proposed in \cite{gadjov2022exact} by introducing an output strictly passive system to a gradient-play scheme.
With similar ideas, the authors of \cite{gao2023second} explored a second-order mirror descent (MD2) dynamics, which converged to a variationally stable state without techniques such as time-averaging or discounting.
However, both \cite{gadjov2018passivity} and \cite{gao2023second} only dealt with NE seeking.

Inspired by the above observations, we consider seeking GNEs over networks for noncooperative games with nonlinear coupled constraints. Taking a typical gradient-play dynamics as a starting point, we focus on analyzing and generalizing the dynamics  under a systematic passivity-based framework.
Our main contributions are summarized as follows.

\begin{itemize}
\item We reveal the convergence of a well-known gradient-play dynamics from a passivity-based perspective, and address that the strict monotonicity assumption on pseudo-gradients plays a important role on its exact convergence. This is an alternative result to \cite{bianchi2021continuous, zeng2019generalized, cherukuri2016asymptotic}, in that the passivity properties are explored.

\item We propose two novel dynamics by introducing parallel feedforward compensators (PFCs) and output strictly passive compensators (OFCs) to the standard gradient-play scheme.
Both of them reach exact GNEs in merely monotone regimes if the PFCs are strictly passive or the OFCs are output strictly passive. They cover some existing methods, including the modified primal-dual algorithms for constrained optimization  in \cite{wang2010control, yamashita2020passivity} and HA for NE seeking in \cite{gadjov2022exact},
and meanwhile, provide ideas on how to design distributed GNE seeking algorithms with convergence guarantees.

\item We explore a unifying framework to generalize the gradient-play dynamics, and develop a class of explicit dynamics for distributed GNE seeking via substituting integrators in the gradient-play dynamics by LTI systems.
We surprisingly find that the dynamics achieves GNEs for strictly monotone games if the LTI systems are passive, and also solve a group of regulator equations.

\item We extend all the proposed dynamics to seek GNEs under partial-decision information settings, and show that
the convergence results still hold under similar assumptions as those of \cite{gadjov2018passivity, tatarenko2020geometric, bianchi2021continuous}.
\end{itemize}

A short version of this work appeared in \cite{li2024passivity},  where a passivity-based gradient-play dynamics was proposed by introducing PFCs.
We should mention that the dynamics proposed in \cite{li2024passivity} is just a special case of (\ref{alg:GGP:forward}).
Different from \cite{li2024passivity}, OFCs are also employed to develop another dynamics in this paper, and a unifying passivity-based framework is explored to generalize the gradient-play scheme.
Moreover, we extend the results to partial-decision information settings.

The rest of this paper is organized as follows.
In Section \ref{sec:background}, we introduce some preliminaries, and formulate the distributed GNE seeking problem.
We revisit a typical gradient-play dynamics based on passivity in Section \ref{sec:GP}.
Then we propose two novel dynamics by introducing PFCs and OFCs to a gradient-play scheme with convergence analysis in Sections \ref{sec:forward} and \ref{sec:back}.
Following that, we explore a unifying framework to generalize the gradient-play dynamics in Section \ref{sec:generalized}.
We extend the theoretical results to partial-decision information settings in Section \ref{sec:partial-info}, and provide illustrative simulations in Section \ref{sec:simulation}.
Finally, we close this paper with concluding remarks in Section \ref{sec:conclusion}.

\emph{Notation:} 
Let $\mathbb R^m$, $\mathbb R_+^m$ and $\mathbb R^{m \times n}$  be the set of $m$-dimensional real column vectors,
$m$-dimensional nonnegative real column vectors, 
and $m$-by-$n$ dimensional real matrices.
Denote $1_m$ ($0_m$) as the $m$-dimensional column vector with all entries of $1$ ($0$), and $I_n$ as the $n$-by-$n$ identity matrix.
We simply write $\mathbf{0}$ for vectors/matrices of zeros with appropriate dimensions when there is no confusion.
The Euclidean inner product of $x$ and $y$ is $x^{\rm T} y$ or $\langle x, y \rangle$.
Let $(\cdot)^{\rm T}$, $\otimes$ and $\Vert \cdot \Vert$ be the transpose, the Kronecker product and the Euclidean norm.
Denote ${\rm Ker}(\cdot)$, ${\rm Im}(\cdot)$ and ${\rm rank}(\cdot)$ as the kernel space, the image space and the rank of a matrix, respectively.
Denote ${\rm dim}(\cdot)$ as the dimension of a linear space.
We call a matrix nonnegative if all its entries are no smaller than zero.
For $X \in \mathbb R^{n \times n}$, $X \succ \mathbf{0}$ means that $X$ is positive definite.
Take $\mathcal I = \{1, \dots N\}$.
For  $x_i \in \mathbb R^{n_i}$,  
we define ${\rm col}\{x_i\}_{i \in \mathcal I} = [x_1^T, \dots, x_N^T]^T  \in \mathbb R^n$, where $n = \sum_{i \in \mathcal I} n_i$.
For $X_i \in \mathbb R^{p_i \times q_i}$, ${\rm blkdiag}\{X_i\} _{i \in \mathcal I}\in \mathbb R^{p \times q}$ is a block diagonal matrix defined by $X_i$, where
$p = \sum_{i \in \mathcal I} p_i$ and $q = \sum_{i \in \mathcal I} q_i$.
Given a differentiable function $J(x, y)$, $\nabla_x J(x, y)$ is the partial-gradient of $J$ with respect to $x$.

\section{BACKGROUND}
\label{sec:background}

In this section, we introduce some necessary concepts, and then, formulate the distributed GNE seeking problem.

\subsection{Mathematical Preliminaries}

Consider a multi-agent network modeled by an undirected graph $\mathcal G(\mathcal I, \mathcal E, \mathcal A)$, where $\mathcal I = \{1, \dots, N\}$ is the node set, $\mathcal E \subset \mathcal I \times \mathcal I$ is the edge set, and $\mathcal A = [a_{ij}] \in \mathbb R^{N \times N}$ is the adjacency matrix such that $a_{ij} = a_{ji} > 0$ if $(i, j) \in \mathcal E$, and otherwise, $a_{ij} = 0$.
Suppose there are no self-loops in $\mathcal G$, i.e., $a_{ii} = 0, \forall i \in \mathcal I$.
The Laplacian matrix $\mathcal L$ is $\mathcal L = \mathcal D - \mathcal A$, where $\mathcal D = {\rm diag}\{d_i\}\in \mathbb R^{N \times N}$, and $d_i = \sum_{j \in \mathcal I} a_{ij}$.
Node $j$ is a neighbor of node $i$ if and only if $(i, j) \in \mathcal E$. 
Let $\mathcal I_i = \{j|(i, j)\in \mathcal E\}$ be the set of node $i$’s neighbors. 
Graph $\mathcal G$ is connected if there exists a path between any pair of distinct nodes.
If $\mathcal G$ is connected, then $\mathcal L = \mathcal L^{\rm T}$, ${\rm rank}(\mathcal L) = N - 1$, and ${\rm Ker}(\mathcal L) = \{k1_N: k \in \mathbb R\}$.


Let $\Omega \subset \mathbb R^m$ be a convex set such that
$\lambda x + (1-\lambda)y \in \Omega, 
\forall x , y \in \Omega, \forall \lambda \in [0, 1]$.
Its tangent cone at $x \in \Omega$ is
$\mathcal T_{\Omega}(x) = \big\{ 
\lim_{k \to \infty} \frac{x_k - x}{\tau_k} ~|~
x_k \in \Omega, x_k \to x, \tau_k > 0, \tau_k \to 0 \big\}$,
while its normal cone is
$\mathcal N_{\Omega}(x) =
\{v \in \mathbb R^m | v^{\rm T}(y-x) \le 0, \forall y \in \Omega\}$.

Take
${\rm proj}_{\Omega}(x) = {\rm argmin}_{y \in \Omega} \Vert y - x\Vert, \forall x \in \mathbb R^m$.
A differentiated projection operator $\Pi_{\Omega}(x, \cdot)$ on $\mathcal T_{\Omega}(x) $ is defined by $\Pi_{\Omega}(x, v) := {\rm prox}_{\mathcal T_{\Omega(x)}}(v) = \lim_{h \to 0^+} \frac {{\rm proj}_{\Omega}(x + h v) - x}{h}$.
It follows from \cite[Theorem 6.30]{bauschke10convex} that
\begin{equation}
\label{proj:decomp}
v = {\rm proj}_{\mathcal T_{\Omega(x)}}(v) + {\rm proj}_{\mathcal N_{\Omega(x)}}(v), ~\forall v \in \mathbb R^m.
\end{equation}

Given a differentiable function $f: \Omega \rightarrow \mathbb R$, $\nabla f(x)$ denotes its gradient at $x$. Then $f$ is convex if $\Omega$ is a convex set and
$f(y) \ge f(x) + \langle  \nabla f(x), y -x\rangle, \forall x, y \in \Omega, \forall \lambda \in [0, 1]$.

An operator (or mapping) $F: \Omega \subset \mathbb R^n \rightarrow \mathbb R^n$ is monotone  if $\langle x - y, F(x) - F(y)\rangle \ge 0, \forall x, y \in \Omega$, strictly monotone if the strict inequality holds for all $x \not= y$, $\mu$-strongly monotone if $\langle x - y, F(x) - F(y)\rangle \ge \mu \Vert x - y \Vert^2, \forall x, y \in \Omega$, and $\nu$-hypomonotone if 
$\langle x - y, F(x) - F(y)\rangle \ge -\nu \Vert x - y \Vert^2, \forall x, y \in \Omega$.
$F$ is $\theta$-Lipschitz continuous if $\Vert F(x) - F(y)\Vert \le \theta \Vert x - y \Vert, \forall x, y \in \Omega$.
Note that if $f: \Omega \rightarrow \mathbb R^n$ is a convex function, then $\nabla f$ is monotone.

Consider a nonlinear system given by
\begin{equation}
\label{Sys:Nonlinear}
\dot x = f(x, u), ~
y = h(x, u)
\end{equation}
where $x \in \mathbb R^n$, $u, y \in \mathbb R^m$, $f: \mathbb R^n \times \mathbb R^m \rightarrow \mathbb R^n$ is locally Lipschitz, $h: \mathbb R^n \times \mathbb R^m \rightarrow \mathbb R^m$ is continuous, $f(\mathbf{0}, \mathbf{0}) = \mathbf{0}$, and $h(\mathbf{0}, \mathbf{0}) = \mathbf{0}$. 
If there exists a continuous differentiable storage function $V(x)$ such that
$\dot V = \nabla V(x)^T f(x, u) \le u^T y, ~\forall (x, u) \in \mathbb R^n \times \mathbb R^m$,
then (\ref{Sys:Nonlinear}) is a passive system.
Moreover, it is passive lossless if $\dot V = u^T y$, strictly passive  if $\dot V \le u^T y - \phi (x)$ for some positive definite function $\phi$,
and output strictly passive if
$\dot V \le u^T y - y^T \psi(y)$, where $y^T \psi(y) > 0, \forall y \not= \mathbf{0}$.
Throughout this paper, whenever we refer to a passive, strictly passive or output strictly passive system, we alway assume the storage function $V$ is positive definite and radially unbounded.

Let $G(s)$ be a $m \times m$ proper rational transfer function matrix with
a minimal (controllable and observable) state-space realization as
\begin{equation}
\label{Sys:LTI}
\dot x = Ax + Bu, ~y = Cx + Du,
\end{equation}
where $A \in \mathbb R^{n \times n}$, $B \in \mathbb R^{n \times m}$, $C \in \mathbb R^{m \times n}$ and $D \in \mathbb R^{m \times m}$.
We call $G(s)$ positive real (PR) if 
\begin{enumerate}[{\textbullet}]
\item poles of all elements of $G(s)$ are in ${\rm Re}[s] \le 0$;

\item for all real $\omega$ for which $\omega$ is not a pole of any element of $G(s)$, the matrix $G(j \omega) + G^T(-j\omega)$ is positive semidefinite;

\item any pure imaginary pole $j \omega$ of any element of $G(s)$ is a simple pole, and the residue matrix $\lim_{s \rightarrow j \omega} (s-j \omega)G(s)$ is positive semidefinite Hermitian.
\end{enumerate}
Furthermore, $G(s)$ is strictly positive real (SPR) if
$G(s - \epsilon)$ is PR for some $\epsilon > 0$.

Suppose $G(s)$ is PR. In light of \cite[Lemmas. 6.2 and 6.4]{khalil2002nonlinear}, the linear time-invariant (LTI) system (\ref{Sys:LTI}) is passive, and there exist real matrices $P \succ \mathbf{0}$, $L$ and $W$ such that
$A^T P + PA = - L^T L$,
$PB = C^T - L^T W$, and
$W W^T = D + D^T$.
Furthermore, if $G(s)$ is SPR, then (\ref{Sys:LTI}) is strictly passive, and there exist constant $\epsilon > 0$,  real matrices $P \succ \mathbf{0}$, $L$ and $W$ such that
$A^T P + PA = - L^T L - \epsilon P$,
$PB = C^T - L^T W$ and
$W W^T = D + D^T$
by \cite[Lemma. 6.3]{khalil2002nonlinear}.

Referring to \cite{kottenstette2014relationships}, we say (\ref{Sys:LTI}) is output strictly passive with respect to a positive definite storage function $V$ if $\dot V \le u^T y - \delta y^T y$ for some $\delta > 0$.

\subsection{Game Setup}

Consider a set $\mathcal I = \{1, \dots, N\}$ of $N$ players (agents) involved in a game. 
Player $i$ controls its action $x^i \in \mathbb R^{n_i}$, where $\sum_{i = 1}^N n_i = n$.
Let $x = (x^i, x^{-i}) \in \mathbb R^n$ be the $N$-tuple of all agents' actions, where $x^{-i}$ is the $(N-1)$-tuple of all agents' actions except agent $i$'s decision. Alternatively,
$x = {\rm col}\{x^i\}_{i \in \mathcal I} \in \mathbb R^n$.
Player $i$ aims at minimizing its local cost function $J_i(x^i, x^{-i}): \mathbb R^n \rightarrow \mathbb R$, which depends both on its local strategy $x^i$ and the actions of its opponents $x^{-i}$.
Furthermore, there exist separable coupling constraints $X = \{x \in \mathbb R^n~|~g(x) \le 0_m\}$, where $g(x) := \sum_{i\in \mathcal I} g_i(x^i)$, and $g_i: \mathbb R^{n_i} \rightarrow \mathbb R^m$ is a private function only known by agent $i$.
Given $x^{-i}$, the feasible decision set of agent $i$ is
$X_i(x^{-i}) = \{x^i \in \mathbb R^{n_i}: (x^i, x^{-i}) \in X\}$.
Player $i$ aims to find the best response strategy 
\begin{equation}
\label{formulation}
\min_{x^i\in \mathbb R^{n_i}} J_i(x^i, x^{-i}),~
{\rm s.t.}~x^i \in X_i(x^{-i}).
\end{equation}

A collective strategy $x^* = (x^{i, *}, x^{-i,*})$ is called a generalized Nash equilibrium (GNE) if
\begin{equation*}
\label{def:GNE}
x^{i, *} \in {\rm argmin}_{x^i} J_i(x^i, x^{-i,*}),
~{\rm s.t.}~(x^i, x^{-i,*}) \in X,~\forall i \in \mathcal I.
\end{equation*}
Particularly, $x^*$ is called an Nash equilibrium (NE) if  there are no coupled constraints $X$, i.e., $g_i(x^i) \equiv \mathbf{0}$.

To ensure the well-posedness of (\ref{formulation}), we make the following standard assumption.

\begin{assumption}
\label{ass:convex}
For every $i \in \mathcal I$, $J_i$ is continuously differentiable and convex in 
$x^i$, given $x^{-i}$, and $g_i$ is continuous, differentiable and convex.
Furthermore, $X$ is non-empty and satisfies the Slater’s constraint qualification.
\end{assumption}

Under Assumption \ref{ass:convex}, $x^*$ is a GNE if and only if 
the following Karush-Kuhn-Tucker (KKT) conditions hold \cite{facchinei2010generalized}:
\begin{equation}
\label{vGNE}
\begin{aligned}
0_{n_i} = & ~\nabla_{x^i} J_i(x^{i, *}, x^{-i,*}) + 
\nabla g_i(x^{i, *})^T \lambda^{i, *}, \\
0_m \in & - g(x^*) + \mathcal N_{\mathbb R^m_+}(\lambda^{i, *}),
\end{aligned}
\end{equation}
where $\lambda^{i, *} \in \mathbb R^m$ is the Lagrangian multiplier of agent $i$.

Given $x^*$ as a GNE of (\ref{formulation}), the corresponding Lagrangian multipliers may be different for the players, i.e., $\lambda^{i,*} \not= \lambda^{j,*}$ for $i \not= j$.
In this work, we focus on seeking a GNE with the same Lagrangian multiplier, named variational GNE (v-GNE), i.e., $\lambda^{i,*} = \lambda^*_c$, $\forall i\in \mathcal I$ \cite{facchinei2010generalized, kulkarni2012variational}, and we simply
call it a GNE.

\begin{remark}
In fact, (\ref{formulation}) is a general model, and appears in applications including coverage maximization, electricity market, and optical networks \cite{hu2022distributed, bianchi2021continuous}. It covers the formulation in \cite{gadjov2022exact, yi2019operator} by allowing  nonlinear coupled constraints.
\end{remark}

\section{GRADIENT-PLAY DYNAMICS}
\label{sec:GP}

In this section, we revisit a well-known distributed gradient-play dynamics to solve (\ref{formulation}), and address its convergence from a passivity-based perspective.

We consider that agent $i$ only knows its local data, i.e., $J_i$ and $g_i$, which contains its own private information, but has knowledge of all its opponents' decisions $x^{-i}$.
This is a full-decision information setting as discussed in 
\cite{hu2022distributed, yi2019operator, zeng2019generalized}, where agent $i$ can compute the exact partial-gradient $\nabla_{x^i} J_i(x^i, x^{-i})$ at each step.

Here we introduce the algorithm notations. Agent $i$ controls its local strategy $x^i \in \mathbb R^{n_i}$, as well as a local copy of  multiplier $\lambda^i \in \mathbb R^m$ to estimate $\lambda_c^*$.
To ensure all local multipliers $\lambda^i$ reaching consensus, we introduce a auxiliary variable $z^i \in \mathbb R^m$ for agent $i$.
Suppose all agents communicate through an undirected graph $\mathcal G_c$. Specifically, agent $i$ can receive  $\{\lambda^j, z^j\}$ from agent $j$ if and only if $j \in \mathcal I_i$, where $\mathcal I_i$ is the neighbor set of agent $i$. We make the following assumption on $\mathcal G_c$.
\begin{assumption}
\label{ass:graph}
The undirected graph $\mathcal G_c$ is connected.
\end{assumption}

Then a typical gradient-play dynamics for (\ref{formulation}) is 
\begin{equation}
\label{alg:GP:full}
\left\{
\begin{aligned}
\dot x^i = & - \nabla_{x^i} J_i(x^i, x^{-i}) - \nabla g_i(x^i)^T \lambda^i, \\ 
\dot z^i = & \sum_{j \in \mathcal I_i} a_{ij}(\lambda^i - \lambda^j),  \\ 
\dot \lambda^i =& \Pi_{\mathbb R^m_+}\big[\lambda^i, g_i(x^i) \!- \! \sum_{j \in \mathcal I_i} a_{ij}(z^i - z^j) \!-\! \sum_{j \in \mathcal I_i} a_{ij}(\lambda^i - \lambda^j)\big]\! \\
\end{aligned}
\right.
\end{equation}
where $a_{ij}$ is the $(i, j)$-th entry of the adjacency matrix of $\mathcal G_c$, $x^i(0) \in \mathbb R^{n_i}$, $z^i(0) \in \mathbb R^{m}$, and $\lambda^i(0) \in \mathbb R^m_+$.

To rewrite (\ref{alg:GP:full}) in a compact form, we define a pseudo-gradient mapping $F$ as
\begin{equation}
\label{pseudo-gradient}
F(x) : = {\rm col}\{\nabla_{x^i} J_i(x^i, x^{-i})\}_{i \in \mathcal I} \in \mathbb R^n.
\end{equation}
Denote by 
$\nabla G(x) := {\rm blkdiag}\{\nabla g_i(x^i)\}_{i \in \mathcal I}  \in \mathbb R^{\tilde m \times n}$,
$G(x) := {\rm col} \{g_i(x^i)\}_{i \in \mathcal I} \in \mathbb R^{\tilde m}$,
$\lambda := {\rm col} \{\lambda^i\}_{i \in \mathcal I} \in \mathbb R^{\tilde m}$, 
$z := {\rm col} \{z^i\}_{i \in \mathcal I} \in \mathbb R^{\tilde m}$,  and
$L := \mathcal L \otimes I_m \in \mathbb R^{\tilde m \times \tilde m}$,
where $\tilde m = Nm$, and $\mathcal L$ is the Laplacian matrix of $\mathcal G_c$. Therefore, dynamics (\ref{alg:GP:full}) reads as 
\begin{equation}
\label{alg:GP}
\left\{
\begin{aligned}
\dot x = & - F(x) - \nabla G(x)^T \lambda, \\
\dot z = & L \lambda, \\
\dot \lambda =& \Pi_{\mathbb R^{\tilde m}_+} \big[\lambda, G(x) - L z - L \lambda \big]
\end{aligned}
\right.
\end{equation}
where $x(0) \in \mathbb R^n$, $z(0) \in \mathbb R^{\tilde m}$,
and $\lambda(0) \in \mathbb R^{\tilde m}_+$.

\begin{remark}
By the viability theorem \cite{aubin1984differential},
$\lambda(t) \in \mathbb R_+^{\tilde m}$ for all $t \ge 0$ due of the projection operator
$\Pi_{\mathbb R^{\tilde m}_+}[\lambda, \cdot]$.
In the absence of the constraints $g(x) \le 0_m$, (\ref{alg:GP}) degenerates into the gradient-play dynamics in \cite{gadjov2018passivity, shamma2005dynamic}.
If $J_i$ only depends on $x^i$,  it is consistent with the primal-dual method for constrained optimization in \cite{cherukuri2016asymptotic}.
In fact, (\ref{alg:GP}) has been proposed for distributed GNE seeking with convergence analysis in \cite{bianchi2021continuous, zeng2019generalized}.
As a comparison, we focus on analyzing and generalizing the dynamics based on the concept of passivity.
\end{remark}

Referring to Theorem $4.1$ in \cite{zeng2019generalized} or Lemma $2$ in \cite{bianchi2021continuous}, the next  lemma addresses the relationship between equilibria of (\ref{alg:GP}) and GNEs of (\ref{formulation}).

\begin{lemma}
\label{lem:GP:KKT-EQ}
Let Assumptions \ref{ass:convex} and \ref{ass:graph} hold, and $F$ be monotone. A profile 
$x^*$ is a GNE of (\ref{formulation}) if and only if there exists $(\lambda^*, z^*)$ such that $(x^*, \lambda^*, z^*)$
is an equilibrium point of (\ref{alg:GP}), i.e.,
\begin{equation}
\label{KKT-EQ}
\left\{
\begin{aligned}
\mathbf{0} =~& F(x^*) + G(x^*)^T \lambda^*, \\
\mathbf{0}=~& L \lambda^* , \\
\mathbf{0} \in~& G(x^*) - L z^* - L \lambda^* - \mathcal N_{\mathbb R^{\tilde m}_+}(\lambda^*).
\end{aligned}
\right.
\end{equation}
\end{lemma}

Next, we put (\ref{alg:GP}) in a block diagram representation, and analyze its properties. Take $v_x := {\rm col}\{v_x^i\}$, $v_z := {\rm col}\{v_z^i\}$ and $v_{\lambda} := {\rm col}\{v_{\lambda}^i\}$,
where $v_x^i = - \nabla_{x^i} J_i(x^i, x^{-i}) - \nabla g_i(x^i)^T \lambda^i$,
$v_z^i = \sum_{j \in \mathcal I_i} a_{ij}(\lambda^i - \lambda^j)$
and $v_{\lambda}^i = g_i(x^i) -  \sum_{j \in \mathcal I_i} a_{ij}(z^i - z^j) - \sum_{j \in \mathcal I_i} a_{ij}(\lambda^i - \lambda^j)$.

\begin{figure}[htp]
\centering
\includegraphics[scale=0.42]{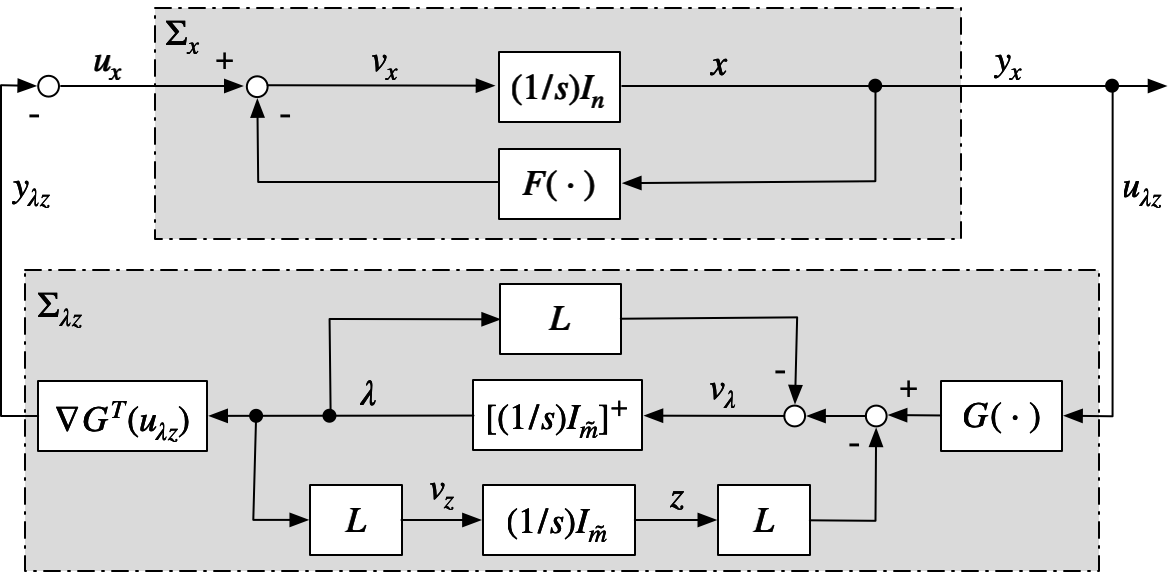}
\caption{Block diagram of dynamics (\ref{alg:GP}).}
\label{fig:GP}
\end{figure}

It is clear that $x$ in (\ref{alg:GP}) can be represented by a bank of integrators, i.e., $x(s) = [(1/s) I_n] v_x(s)$. Similarly,  $z(s) = [(1/s) I_{\tilde m}] v_z(s)$.
To handle $\lambda$ obtained via the projection operator $\Pi_{\mathbb R_+^{\tilde m}}[\lambda, \cdot]$ in cascade with a bank of integrators, we introduce the notation $\lambda(s) = [(1/s) I_{\tilde m}]^+ v_{\lambda}(s)$. Consequently, (\ref{alg:GP}) can be represented by the block diagram shown in  Fig. \ref{fig:GP}, which inspires us to decompose (\ref{alg:GP}) into two interconnected subsystems  $\Sigma_x$ and $\Sigma_{\lambda z}$,
where
\begin{equation}
\label{sigma:x}
\Sigma_x: ~
\dot x = - F(x) + u_x, 
~y_x = x
\end{equation}
and moreover,
\begin{equation}
\label{sigma:lz}
\Sigma_{\lambda z}: 
\left \{
\begin{aligned}
\dot \lambda & = \Pi_{\mathbb R_+^{\tilde m}}\big[\lambda, 
G(u_{\lambda z}) - L z - L \lambda \big], \\
\dot z& =L \lambda, \\
y_{\lambda z} &= \nabla G(u_{\lambda z})^T \lambda.
\end{aligned}
\right.
\end{equation}	

Let $(x^*, \lambda^*, z^*)$ be an equilibrium point of (\ref{alg:GP}), and
\begin{equation}
\begin{aligned}
\label{def:notation}
&\tilde u_x := u_x - u^*_x,
&\tilde u_{\lambda z} &:= u_{\lambda z} - u_{\lambda z}^*,\\
&\tilde y_x := y_x - y_x^*, 
&\tilde y_{\lambda z} &:= y_{\lambda z} - y_{\lambda z}^*,
\end{aligned}
\end{equation}
where $u_{\lambda z}^* = y_x^* = x^*$ and $u_{x}^* = - y_{\lambda z}^* = - \nabla G(x^*)^T \lambda^*$.

Convergence of (\ref{alg:GP}) is established in the following theorem, whose proof can be found in \cite{li2024passivity}.

\begin{theorem}
\label{thm:GP:convergence}
Consider dynamics (\ref{alg:GP}).
Let Assumptions \ref{ass:convex} and \ref{ass:graph} hold, and $F$ be monotone. 
\begin{enumerate}[(i)]
\item System $\Sigma_x$ (\ref{sigma:x}) is passive from $\tilde u_x$ to $\tilde y_x$ with respect to the storage function $V_x = \frac 12 \Vert x - x^* \Vert^2$.
\item System $\Sigma_{\lambda z}$ (\ref{sigma:lz}) is passive from $\tilde u_{\lambda z}$ to $\tilde y_{\lambda z}$ with respect to the storage function 
$V_{\lambda z} = \frac 12 \Vert  \lambda - \lambda^* \Vert^2 +  \frac 12 \Vert  z - z^* \Vert^2.$

\item If $F$ is strictly monotone, then
the trajectory of $\big(x(t), \\ \lambda(t), z(t)\big)$ converges to an equilibrium point of (\ref{alg:GP}), where  the $x(t)$ component approaches the GNE of (\ref{formulation}).
\end{enumerate}
\end{theorem}

\begin{remark}
\label{rmk: GP}
Theorem \ref{thm:GP:convergence} is an alternative result to  \cite{bianchi2021continuous, zeng2019generalized, cherukuri2016asymptotic}, in that passivity properties of $\Sigma_x$ and $\Sigma_{\lambda z}$ are analyzed.
Note that the strict monotonicity assumption on $F$ plays a critical role on the asymptotic convergence of (\ref{alg:GP}). If the assumption fails, (\ref{alg:GP}) may be passive lossless, and cannot be convergent. 
The next example is given for illustration.
\end{remark}

\begin{example}
\label{Ex1}
Consider a two-players zero-sum game problem without constraints, where $J_1(x^1, x^2) = x^1 x^2$ and $J_2(x^2, x^1) = - x^1 x^2$.
Then (\ref{alg:GP}) is cast into
\begin{equation*}
\dot x^1 = - x^2,~\dot x^2 = x^1. 
\end{equation*}
Clearly, the NE $x^*$ is $[0, 0]^T$, but $x(t)$ will cycle around the NE and never converge if $x(0) \not= x^*$.
\end{example}

From Theorem \ref{thm:GP:convergence} and Remark \ref{rmk: GP}, we conclude that passivity is a  powerful tool to prove the convergence of (\ref{alg:GP}), but also provides us with a better understanding of the dynamics by explaining why it can work under certain settings.
It is natural to consider how to generalize 
(\ref{alg:GP}) in a systematic passivity-based framework.
Throughout this paper, we are interested in the following questions.
\begin{itemize}
\item Can we develop any novel dynamics based on (\ref{alg:GP}) for merely monotone games?

\item Can we explore a unifying framework, and design a class of explicit dynamics to generalize (\ref{alg:GP}) resorting to passivity?

\item How to extend our proposed dynamics and theoretical results to partial-decision information settings?
\end{itemize}

\section{PARALLEL FEEDFORWARD COMPENSATION}
\label{sec:forward}

In this section, we propose a novel gradient-play dynamics by adding parallel feedforward compensators (PFCs) to (\ref{alg:GP}). Then we show that the dynamics can achieve exact convergence in merely monotone regimes.
Furthermore, we relate it with two existing algorithms.

\subsection{Algorithm Design}

Referring to \cite[Theorem. 2.10]{sepulchre2012constructive}, 
the passivity of a system will be preserved if a passive compensator is added in parallel. 
The parallel feedforward compensator (PFC) may also adjust the performance of the original system.
Recent applications of PFCs can be found in \cite{kim2016design, yamashita2020passivity, li2022parallel}.
Inspired by these observations, we consider enhancing the convergence of (\ref{alg:GP}) by designing suitable PFCs for $\Sigma_x$ (\ref{sigma:x}) and $\Sigma_{\lambda z}$ (\ref{sigma:lz}).

As shown in Fig. \ref{fig:GP}, only passive integrators $(1/s) I_{n_i}$ are involved in the evolution of $x^i$.
For each $i \in \mathcal I$, we introduce a PFC to $(1/s) I_{n_i}$.
The PFC is a LTI system $H_i^f(s)$ with a minimal state-space realization as
\begin{equation}
\label{PFC:x}
\dot \tau^i_x = \Phi_i^f \tau^i_x + \Theta_i^f v^i_x, ~
\eta^i_x  = \Psi_i^f  \tau^i_x,
\end{equation}
where $\Phi_i^f  \in \mathbb R^{p_i^f \times p_i^f}$ is Hurwitz,
$\Theta_i^f \in \mathbb R^{p_i^f \times n_i}$ is full-column rank,
and $\Psi_i^f  \in \mathbb R^{n_i \times p_i^f}$ is full-row rank.
Let $\tau_x := {\rm col}\{\tau_x^i\}_{i \in \mathcal I}$,
$\eta_x := {\rm col}\{\eta_x^i\}_{i \in \mathcal I}$,
$\Phi^f := {\rm blkdiag}\{\Phi_i^f\}_{i \in \mathcal I}$,
$\Theta^f := {\rm blkdiag}\{\Theta_i^f\}_{i \in \mathcal I}$, and
$\Psi^f := {\rm blkdiag}\{\Psi_i^f\}_{i \in \mathcal I}$.
Then a state-space representation for $x$ is
\begin{equation}
\label{forward:x:statespace}
\left \{
\begin{aligned}
\dot \rho_x& = v_x, \\
\dot \tau_x &= \Phi^f \tau_x + \Theta^f v_x, ~
\eta_x = \Psi^f  \tau_x, \\
x   &= \rho_x + \eta_x.
\end{aligned}
\right.
\end{equation}

Similarly, for $z^i$ in Fig. \ref{fig:GP}, we add a PFC 
$\hat H_i^f(s)$ to $(1/s) I_m$. A minimal state-space realization of $\hat H_i^f(s)$ is
\begin{equation}
\label{PFC:z}
\dot \tau^i_z = \hat \Phi_i^f \tau^i_z + \hat \Theta_i^f v^i_z, ~
\eta^i_z  = \hat \Psi_i^f  \tau^i_z,
\end{equation}
where $\hat \Phi_i^f  \in \mathbb R^{q_i^f \times q_i^f}$ is Hurwitz,
$\hat \Theta_i^f \in \mathbb R^{q_i^f \times m}$ is full-column rank, and $\hat \Psi_i^f  \in \mathbb R^{m \times q_i^f}$ is full-row rank.
Let $\tau_z := {\rm col}\{\tau_z^i\}_{i \in \mathcal I}$,
$\eta_z := {\rm col}\{\eta_z^i\}_{i \in \mathcal I}$,
$\hat \Phi^f := {\rm blkdiag}\{\hat \Phi_i^f\}_{i \in \mathcal I}$,
$\hat \Theta^f := {\rm blkdiag}\{\hat \Theta_i^f\}_{i \in \mathcal I}$, and
$\hat \Psi^f := {\rm blkdiag}\{\hat \Psi_i^f\}_{i \in \mathcal I}$.
Consequently, a state-space representation for $z$ is
\begin{equation}
\label{forward:z:statespace}
\left \{
\begin{aligned}
\dot \rho_z& = v_z, \\
\dot \tau_z &=\hat \Phi^f \tau_z + \hat \Theta^f v_z, ~
\eta_z = \hat \Psi^f  \tau_z,\\
z & = \rho_z + \eta_z.
\end{aligned}
\right.
\end{equation}

Compared with $x^i$ and $z^i$, it is more challenging to handle $\lambda^i$ due to the projection operator $\Pi_{\mathbb R_+^{m}}[\lambda^i, \cdot]$.
For convenience, we define a similar notation as $[(1/s) I_m]^+$ in Fig. \ref{fig:GP}.

Let $G(s)$ be the transfer function of a minimal LTI system 
$\dot x = Ax + Bu, ~y = Cx + Du$,
where $u, y \in \mathbb R^m$, and $x \in \mathbb R^n$. We denote by
$[G(s)]^+$ to mean $y(s) = G(s) u(s)$ under the constraint $y(t) \ge \mathbf{0}$ for all $t \ge 0$. Specifically, in time domain, a state-space realization of $[G(s)]^+$ is 
\begin{equation}
\label{def:proj:tranf}
\dot x = \Pi_{\mathbb R^n_+}[x, Ax + Bu], ~
y ={\rm max}\{\mathbf{0}, C x + Du\}.
\end{equation}

For $\lambda^i$  in Fig. \ref{fig:GP}, we introduce a PFC 
$[\bar H_i(s)]^+$ for $[(1/s) I_m]^+$. A state-space realization of $[\bar H_i(s)]^+$ is 
\begin{equation}
\label{PFC:l}
\dot \tau^i_{\lambda} =\Pi_{\mathbb R^m_+}[\tau^i_{\lambda} , \bar \Phi^f_i \tau^i_{\lambda}  + \bar \Theta^f_i v^i_{\lambda}] , ~
\eta^i_{\lambda} = {\rm max}\{\mathbf{0}, \bar \Psi^f_i  \tau^i_{\lambda}\}
\end{equation}
where $\bar \Phi_i^f  \in \mathbb R^{r_i^f \times r_i^f}$ is negative definite,
$\bar \Theta_i^f \in \mathbb R^{r_i^f \times m}$ is nonnegative and full-column rank, and $\bar \Psi_i^f = \bar \Theta^{f, T}_i \in \mathbb R^{m \times r_i^f}$.
Let $\tau_{\lambda} := {\rm col}\{\tau_{\lambda}^i\}_{i \in \mathcal I}$,
$\eta_{\lambda} := {\rm col}\{\eta_{\lambda}^i\}_{i \in \mathcal I}$,
$\bar \Phi^f := {\rm blkdiag}\{\bar \Phi_i^f\}_{i \in \mathcal I}$,
$\bar \Theta^f := {\rm blkdiag}\{\bar \Theta_i^f\}_{i \in \mathcal I}$, 
$\bar \Psi^f := {\rm blkdiag}\{\bar \Psi_i^f\}_{i \in \mathcal I}$,
and $r^f = \sum_{i \in \mathcal I} r_i^f$.
Therefore, a state-space representation for $\lambda$ is
\begin{equation}
\label{forward:l:statespace}
\left \{
\begin{aligned}
\dot \rho_{\lambda} &= \Pi_{\mathbb R_+^{\tilde m}}[\rho_{\lambda}, v_{\lambda}], \\
\dot \tau_{\lambda} &=\Pi_{\mathbb R_+^{r^f}}[\tau_{\lambda}, \bar \Phi^f \tau_{\lambda} + \bar \Theta^f v_{\lambda}], ~
\eta_{\lambda} =\max\{\mathbf{0}, \bar \Psi^f  \tau_{\lambda}\}, \\
\lambda &= \rho_{\lambda} + \eta_{\lambda}.
\end{aligned}
\right.
\end{equation}

We make the following assumption on the introduced PFCs.
\begin{assumption}
\label{ass:forward}
For each $i \in \mathcal I$, $H_i^f(s)$ and $\hat H_i^f(s)$ are SPR, i.e., systems (\ref{PFC:x}) and (\ref{PFC:z}) are strictly passive with respect to storage functions $S^f_{\tau_x^i} = \frac 12 \tau_x^{i,T} P_{\tau_x^i} \tau_x^i$, and
$S^f_{\tau_z^i} = \frac 12 \tau_z^{i,T} P_{\tau_z^i} \tau_z^i$,
where $P_{\tau_x^i} \succ \mathbf{0}$ and $P_{\tau_z^i} \succ \mathbf{0}$.
System (\ref{PFC:l}) is strictly passive with respect to the storage function $S^f_{\tau_{\lambda}^i} = \frac 12 \Vert \tau_{\lambda}^i \Vert^2$.
\end{assumption}

\begin{remark}
The strictly passive realness of $H_i^f(s)$ implies that $\Phi^f_i$ in (\ref{PFC:x}) is Hurwitz. For (\ref{PFC:x}), none of its inputs and outputs are redundant since $\Theta^f_i$ and $\Psi^f_i$ are full rank. 
Suppose $\bar \Phi_i^f$ is negative definite, $\bar \Theta_i^f$ is nonnegative and full-column rank, and $\bar \Psi_i^f = \bar \Theta^{f, T}_i$. Then it is straightforward to verify that (\ref{PFC:l}) is strictly passive with respect to $S^f_{\tau_{\lambda}^i}$.
In practice, Assumption \ref{ass:forward} is easy to hold.
Examples of SPR systems can be found in 
\cite[Chapter. 6]{khalil2002nonlinear} \cite{kottenstette2014relationships}, while an instance of (\ref{PFC:l}) is 
$\bar \Phi_i^f  = {\rm diag}\{\bar a_i\} \in \mathbb R^{m \times m}$,
$\bar \Theta_i^f = {\rm diag}\{\bar b_i\} \in \mathbb R^{m \times m}$, and $\bar \Psi^f_i = \bar \Theta^{f, T}_i$,
where $\bar a_i = [-\bar a_{i1}, \dots, - \bar a_{im}]^T \in \mathbb R^m$,
$\bar b_i = [\bar b_{i1}, \dots, \bar b_{im}]^T \in \mathbb R^m$ and $\bar a_{ik}, \bar b_{i k} > 0$ for all $k \in \{1, \dots, m\}$.
\end{remark}

Combining (\ref{forward:x:statespace}), (\ref{forward:z:statespace}) with (\ref{forward:l:statespace}), a novel passivity-based gradient-play dynamics is designed as
\begin{equation}
\label{alg:GGP:forward}
\left \{
\begin{aligned}
\dot \rho_x& = - F(x) - \nabla G(x)^T \lambda, \\
\dot \tau_x &= \Phi^f \tau_x - \Theta^f \big[F(x) + \nabla G(x)^T \lambda\big], ~
\eta_x = \Psi^f  \tau_x \\
x & = \rho_x + \eta_x, \\
\dot \rho_{\lambda} & = \Pi_{\mathbb R_+^{\tilde m}}\big[\rho_{\lambda}, 
G(x) - L z - L \lambda \big], \\
\dot \tau_{\lambda} &=\Pi_{\mathbb R^{r^f}_+}\big[\tau_{\lambda} , \bar \Phi^f \tau_{\lambda} + \bar \Theta^f \big(G(x) - L z - L \lambda \big)\big],\\
\eta_{\lambda} &= {\rm max}\big\{ \mathbf{0}, \bar \Psi^f  \tau_{\lambda} \big\},\\
\lambda &= \rho_{\lambda} + \eta_{\lambda},\\
\dot \rho_z& = L \lambda, ~
\dot \tau_z = \hat \Phi^f \tau_z + \hat  \Theta^f  L \lambda, ~
\eta_z= \hat \Psi^f \tau_z, \\
z & = \rho_z + \eta_z
\end{aligned}
\right.
\end{equation}
where $\rho_x(0) \!\in \! \mathbb R^{n}$, 
$\tau_x(0) \!\in\! \mathbb R^{p^f}$,
$\rho_{\lambda}(0)\! \in\! \mathbb R_+^{\tilde m}$, 
$\tau_{\lambda}(0) \!\in \!\mathbb R_+^{r^f}$
$\rho_z(0) \!\in\! \mathbb R^{\tilde m}$, 
$\tau_z(0) \!\in\! \mathbb R^{q^f}$,
$p^f \!=\! \sum_{i \in \mathcal I} p_i^f$, and
$q^f \!=\! \sum_{i \in \mathcal I} q_i^f$.

Similar to (\ref{alg:GP}), only the first-order information, such as $F$ and $\nabla G$, is used in dynamics (\ref{alg:GGP:forward}). 
It is also clear that (\ref{alg:GGP:forward}) is implemented in a distributed manner.
The next lemma relates the equilibria of (\ref{alg:GGP:forward}) with GNEs of (\ref{formulation}).
\begin{lemma}
\label{lem:forward:KKT-EQ}
Suppose that Assumptions \ref{ass:convex}, \ref{ass:graph} and \ref{ass:forward} hold, and moreover, $F$ is monotone.
\begin{enumerate}[(i)]
\item Let $\big(\rho_x^*, \tau_x^*, \rho_{\lambda}^*, \tau_{\lambda}^*, \rho_z^*, \tau_z^*\big)$ be an equilibrium point of (\ref{alg:GGP:forward}), and $(x^*, \lambda^*, z^*)$ be the corresponding output.
Then ${\rm col}\{\tau_x^*, \tau_{\lambda}^*, \tau_z^*\} = \mathbf{0}$,
${\rm col}\{x^*, \lambda^*, z^*\} = {\rm col}\{\rho_x^*, \rho_{\lambda}^*, \rho_z^*\}$, and the tuple
$(x^*, \lambda^*, z^*)$  satisfies the KKT conditions in (\ref{KKT-EQ}), where $x^*$ is a GNE of (\ref{formulation}).

\item Dynamics (\ref{alg:GGP:forward}) admits at least one equilibrium point $\big(\rho_x^*, \tau_x^*, \rho_{\lambda}^*, \tau_{\lambda}^*, \rho_z^*, \tau_z^*\big)$ such that $\rho_x^* = x^*$, where $x^*$ is a GNE of (\ref{formulation}).
\end{enumerate} 
\end{lemma}

\emph{Proof:}
(i) In light of $\dot \rho_x = \mathbf{0}$, we obtain
$F(x^*) + \nabla G(x^*)^T \lambda^* = \mathbf{0}.$
Recalling $\dot \tau_x = \mathbf{0}$ gives
$ \Phi^f \tau_x^* = \mathbf{0}$.
Since $\Phi^f$ is Hurwitz, $\tau_x^* = \mathbf{0}$, and hence, $x^* = \rho_x^*$.
Similarly, we have $L \lambda^* = \mathbf{0}$, 
$\tau_z^* = \mathbf{0}$,
and $z^* = \rho_z^*$.

Since $\dot \rho_{\lambda} = \mathbf{0}$,
there exists $\zeta_1^*$ such that
$\zeta_1^* = G(x^*) - L z^* - L \lambda^* \in \mathcal N_{\mathbb R^{\tilde m}_+}(\rho_{\lambda}^*)$.
By $\dot \tau_{\lambda} = \mathbf{0}$,
there is $\zeta_2^*$ such that
\begin{equation}
\label{pf:lem:forward:eq:2}
\zeta_2^* = \bar \Phi^f \tau_{\lambda}^* + \bar \Theta^f \zeta_1^* \in \mathcal N_{\mathbb R^{r^f}_+} (\tau_{\lambda}^*).
\end{equation}
Multiplying $\tau_{\lambda}^{*T}$ to both sides of 
(\ref{pf:lem:forward:eq:2}), we derive
\begin{equation}
\label{pf:lem:forward:eq:3}
\tau_{\lambda}^{*T}\zeta_2^* =  \tau_{\lambda}^{*T} \bar \Phi^f \tau_{\lambda}^* +  \tau_{\lambda}^{*T} \bar \Theta^f \zeta_1^*.
\end{equation}
Since $\bar \Phi^f$ is negative definite, $\tau_{\lambda}^{*T} \bar \Phi^f \tau_{\lambda}^* \le 0$, and the equality holds only if $\tau_{\lambda}^* = \mathbf{0}$. Note that $\tau_{\lambda}^{*T}\zeta_2^* = 0$, $\tau_{\lambda}^* \ge \mathbf{0}$, and $\zeta_1^* \le \mathbf{0}$. 
Consequently, (\ref{pf:lem:forward:eq:3}) holds only if 
$\tau_{\lambda}^* = \mathbf{0}$ because $\bar \Theta^f$ is nonnegative. Therefore, $\tau_{\lambda}^* = \mathbf{0}$,
$\lambda^* = \rho_{\lambda}^*$, and moreover,
$G(x^*) - L z^* - L \lambda^* \in \mathcal N_{\mathbb R^{\tilde m}_+}(\lambda^*).$

To sum up, ${\rm col}\{\tau_x^*, \tau_{\lambda}^*, \tau_z^*\} = \mathbf{0}$,
${\rm col}\{x^*, \lambda^*, z^*\} = {\rm col}\{\rho_x^*, \rho_{\lambda}^*, \rho_z^*\}$, and 
$(x^*, \lambda^*, z^*)$  satisfies
(\ref{KKT-EQ}). By Lemma \ref{lem:GP:KKT-EQ}, $x^*$ is a GNE of (\ref{formulation}), and thus, part (i) holds.

(ii) Let $x^*$ be a GNE of (\ref{formulation}).  By Lemma \ref{lem:GP:KKT-EQ}, there exists $(\lambda^*, z^*)$ such that $(x^*, \lambda^*, z^*)$ satisfies (\ref{KKT-EQ}).
Take ${\rm col}\{\tau_x^*, \tau_{\lambda}^*, \tau_z^*\} = \mathbf{0}$. Then ${\rm col}\{\dot \rho_x, \dot \tau_x, \dot \rho_{\lambda}, \dot \rho_z, \dot \tau_z\} = \mathbf{0}$,
and ${\rm col}\{\rho_x^*, \rho_{\lambda}^*, \rho_z^*\} = {\rm col}\{x^*, \lambda^*, z^*\}$.
Note that $\zeta^*_3 := G(x^*) - L z^* - L \lambda^* \le 0$. 
As a result, $\bar \Phi^f \tau_{\lambda}^* + \bar \Theta^f \zeta^*_3 \le 0$ since $\bar \Theta^f$ is nonnegative,
and $\bar \Phi^f \tau_{\lambda}^* + \bar \Theta^f \zeta^*_3 \in \mathcal N_{\mathbb R^{r^f}}(\tau_{\lambda}^*)$, i.e., $\dot \tau_{\lambda} = \mathbf{0}$.
Thus, part (ii) holds. This completes the proof.
$\hfill\square$

\begin{remark}
\label{rmk:forward:KKT}
Lemma \ref{lem:forward:KKT-EQ} indicates that the introduced PFCs do not change the equilibria of (\ref{alg:GP}).
Furthermore, $x(t)$ would reach a GNE of (\ref{formulation}) if the states of (\ref{alg:GGP:forward}) converged to an equilibrium point.
\end{remark}

\subsection{Convergence Analysis}

In this subsection, we analyze the convergence of (\ref{alg:GGP:forward}).

\begin{figure}[htp]
\centering
\includegraphics[scale=0.4]{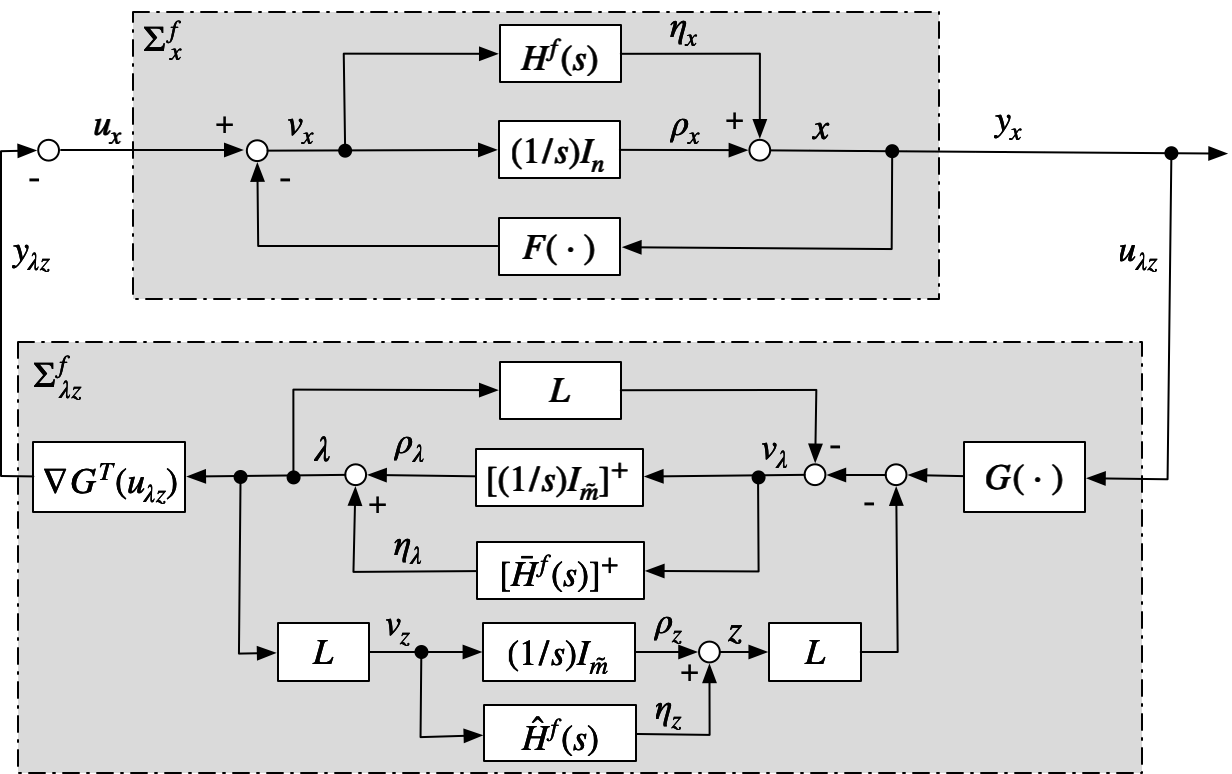}
\caption{Block diagram of dynamics (\ref{alg:GGP:forward}).}
\label{fig:GGP:forward}
\end{figure}

Let $H^f(s) := {\rm blkdiag}\{H^f_i(s)\}_{i \in \mathcal I}$,
$[\bar H^f(s)]^+ := {\rm blkdiag}\{[\bar H^f_i(s)]^+\}_{i \in \mathcal I}$, and
$\hat H^f(s) := {\rm blkdiag}\{\hat H^f_i(s)\}_{i \in \mathcal I}$.
Fig. \ref{fig:GGP:forward} shows the block diagram of dynamics (\ref{alg:GGP:forward}). 
Similar to (\ref{alg:GP}), we decompose 
(\ref{alg:GGP:forward}) into two interconnected subsystems $\Sigma_x^f$ and $\Sigma_{\lambda z}^f$, where
\begin{equation}
\label{sigma:x:forward}
\Sigma_x^f:
\left \{
\begin{aligned}
\dot \rho_x& =  - F(x) + u_x, \\
\dot \tau_x &= \Phi^f \tau_x + \Theta^f \big(- F(x) + u_x \big), ~
\eta_x = \Psi^f  \tau_x,\\
x & = \rho_x + \eta_x, ~
y_x = x
\end{aligned}
\right.
\end{equation}
and moreover,
\begin{equation}
\label{sigma:lz:forward}
\Sigma_{\lambda z}^f \! :\!
\left \{
\begin{aligned}
\dot \rho_{\lambda}& = \Pi_{\mathbb R_+^{\tilde m}}\big[\rho_{\lambda}, G(u_{\lambda z}) - L z - L \lambda \big], \\
\dot \tau_{\lambda} &=\Pi_{\mathbb R_+^{r^f}} \big[\tau_{\lambda}, \bar \Phi^f \tau_{\lambda} + \bar \Theta^f \big(G(u_{\lambda z}) - L z - L \lambda  \big) \big], \\
\eta_{\lambda} &=\max\{\mathbf{0}, \bar \Psi^f  \tau_{\lambda}\},~
\lambda  = \rho_{\lambda} + \eta_{\lambda},\\
\dot \rho_z& = L \lambda, ~
\dot \tau_z =\bar \Phi^f \tau_z + \bar \Theta^f L \lambda, 
~\eta_z = \bar \Psi^f  \tau_z,\\
z & = \rho_z + \eta_z, ~
y_{\lambda z} = \nabla G(u_{\lambda z})^T \lambda.
\end{aligned}
\right.
\end{equation}

Let $\big(\rho(t), \tau(t)\big)$ be a state trajectory of (\ref{alg:GGP:forward}), and $(x(t), \lambda(t), z(t))$ be the corresponding  output trajectory, where $\rho(t) = {\rm col}\{\rho_x(t), \rho_{\lambda}(t), \rho_z(t)\}$ and $\tau(t) = {\rm col}\{\tau_x(t),  \tau_{\lambda}(t), \tau_z(t)\}$.
Suppose $\big(\rho^*, \tau^*\big)$ is an equilibrium point of (\ref{alg:GGP:forward}), and $(x^*, \lambda^*, z^*)$ is the output. By Lemma \ref{lem:forward:KKT-EQ},
$\tau^* = \mathbf{0}$,
$\rho^* = {\rm col}\{x^*, \lambda^*, z^*\}$, and
$(x^*, \lambda^*, z^*)$ satisfies (\ref{KKT-EQ}).
For $\Sigma^f_x$ and $\Sigma_{\lambda z}^f$, we define two storage functions as
\begin{equation}
\label{pf:forward:storage:x}
S^f_x = \frac 12 \Vert \rho_x - x^*\Vert^2 + \frac 12 \tau_x^T P_{\tau_x} \tau_x
\end{equation}
and 
\begin{equation}
\label{pf:forward:storage:lz}
S^f_{\lambda z} = \frac 12 \Vert \rho_{\lambda} - \lambda^* \Vert + \frac 12 \Vert \rho_z - z^*\Vert^2 
+ \frac 12 \Vert \tau_{\lambda} \Vert^2
+ \frac 12 \tau_z^T P_{\tau_z} \tau_z
\end{equation}
where $P_{\tau_x} = {\rm blkdiag}\{P_{\tau^i_x}\}$, 
$P_{\tau_z} = {\rm blkdiag}\{P_{\tau^i_z}\}$
and moreover, $P_{\tau^i_x}$ and $P_{\tau^i_z}$ are given in Assumption \ref{ass:forward}.
We also define $\tilde u_x$, $\tilde y_x$, $\tilde u_{\lambda z}$ and $\tilde y_{\lambda z}$ by (\ref{def:notation}).

The following theorem shows the convergence of  (\ref{alg:GGP:forward}).
\begin{theorem}
\label{thm:forward:convergence}
Consider dynamics (\ref{alg:GGP:forward}). Let Assumptions \ref{ass:convex}, \ref{ass:graph} and \ref{ass:forward} hold, and $F$ be monotone.
\begin{enumerate}[(i)]
\item System $\Sigma^f_x$ (\ref{sigma:x:forward}) is passive from $\tilde u_x$ to $\tilde y_x$ with respect to the storage function $S^f_x$ (\ref{pf:forward:storage:x}).

\item System $\Sigma^f_{\lambda z}$ (\ref{sigma:lz:forward}) is passive from $\tilde u_{\lambda z}$ to $\tilde y_{\lambda z}$ with respect to the storage function $S^f_{\lambda z}$ (\ref{pf:forward:storage:lz}).

\item The state trajectory of $\big(\rho(t), \tau(t) \big)$ converges to an equilibrium point of (\ref{alg:GGP:forward}), where the output trajectory of $x(t)$ approaches a GNE of (\ref{formulation}).
\end{enumerate}
\end{theorem}

\emph{Proof:}
(i) From (\ref{sigma:x:forward}), it follows that
\begin{equation*}
\begin{aligned}
\dot S^f_x &= \big\langle \rho_x - x^*, u_x - F(x) \big \rangle
+ \frac 1 2 \tau_x^T(\Phi^{f,T} P_{\tau_x} + P_{\tau_x} \Phi^f)\tau_x \\
&+ \big\langle \tau_x, P_{\tau_x} \Theta^f \big(u_x - F(x) \big)\big\rangle
\end{aligned}
\end{equation*}
Since $H_i^f$ is SPR, 
$P_{\tau_x} \Theta^f = \Psi^{f,T}$, and
there exists $\delta^f_{\tau_x} > 0$ such that
$\frac 1 2 \tau_x^T(\Phi^{f,T} P_{\tau_x} + P_{\tau_x} \Phi^f)\tau_x \le - \delta^f_{\tau_x}\Vert \tau_x\Vert^2$.
Hence,
$\dot S^f_x \le - \delta^f_{\tau_x}\Vert \tau_x\Vert^2 + 
\big\langle \rho_x + \eta_x - x^*, u_x - F(x) \big\rangle.$
Note that $u_x^* = F(x^*)$ and $y_x = x$. Then
\begin{equation}
\begin{aligned}
\label{pf:forward:conv:1}
\dot S^f_x \le&  - \delta^f_{\tau_x} \Vert \tau_x\Vert^2
-\big\langle x - x^*, F(x) - F(x^*)\big \rangle \\
&+ \big\langle u_x - u_x^*, y_x - y_x^*\big\rangle.
\end{aligned}
\end{equation}

Due to the monotonicity of $F$,
$\dot S^f_x \le \langle \tilde u_x, \tilde y_x \rangle$, and as a result, part (i) holds.

(ii) Recalling (\ref{sigma:lz:forward}) gives
\begin{equation*}
\begin{aligned}
\dot S^f_{\lambda z} =& \big\langle \rho_z - z^*, L \lambda \big\rangle + \frac 12 \tau_z^T \big(\hat \Phi^{f, T} P_{\tau_z} + P_{\tau_z} \hat \Phi^f \big)\tau_z \\
+ & \big\langle \tau_z, P_{\tau_z} \hat \Theta^f L \lambda \big\rangle 
+ \big\langle \tau_{\lambda},  \dot \tau_{\lambda}\big\rangle 
+ \big\langle \rho_{\lambda} - \lambda^*, \dot \rho_{\lambda} \big\rangle.
\end{aligned}
\end{equation*}
Because $\hat H_i^f$ is SPR, 
$P_{\tau_z} \hat\Theta^f = \hat \Psi^{f,T}$, and
there exists $\delta^f_{\tau_z} > 0$ such that
$\frac 1 2 \tau_x^T(\hat \Phi^{f,T} P_{\tau_z} + P_{\tau_z} \hat\Phi^f)\tau_z \le - \delta^f_{\tau_z}\Vert \tau_z\Vert^2$. Thus,
\begin{equation}
\label{pf:forward:conv:2}
\dot S^f_{\lambda z} \le - \delta^f_{\tau_z}\Vert \tau_z\Vert^2 + \big\langle z - z^*, L \lambda\big\rangle
+ \big\langle \rho_{\lambda} - \lambda^*, \dot \rho_{\lambda} \big\rangle + \big\langle \tau_{\lambda},  \dot \tau_{\lambda}\big\rangle .
\end{equation}

According to (\ref{proj:decomp}), we have
\begin{equation}
\begin{aligned}
\label{pf:forward:conv:3}
\langle \rho_{\lambda} - \lambda^*, \dot \rho_{ \lambda} \rangle = \big\langle \rho_{\lambda} - \lambda^*, \big[ G(u_{\lambda z}) - L z - 
L \lambda\big]& \\
- {\rm proj}_{{\mathcal N}_{\mathbb R_+^{\tilde m}}(\rho_{\lambda})}
\big[G(u_{\lambda z}) - L z - 
L \lambda \big] \big\rangle&\\
\le \big\langle \rho_{\lambda} - \lambda^*, 
G(u_{\lambda z}) - L z - L \lambda\big\rangle&.
\end{aligned}
\end{equation}
Similarly,
$\big\langle \tau_{\lambda},  \dot \tau_{\lambda}\big\rangle \le \big\langle \tau_{\lambda}, \bar \Phi^f \tau_{\lambda} + \bar \Theta^f \big(G(u_{\lambda z})-L z-L \lambda \big) \big\rangle$.
Note that $\tau_{\lambda} \ge \mathbf{0}$, $\bar \Psi^f = \bar \Theta^{f, T}$, and  $\bar \Theta^f$ is nonnegative.
As a result, $\eta_{\lambda} = \bar \Psi^f \tau_{\lambda} = \bar \Theta^{f, T} \tau_{\lambda}$.
Since $\bar \Phi^f$ is negative definite, there is $\delta^f_{\tau_{\lambda}} > 0$ such that
\begin{equation}
\label{pf:forward:conv:4}
\big\langle \tau_{\lambda},  \dot \tau_{\lambda}\big\rangle \le -\delta^f_{\tau_{\lambda}}\Vert \tau_{\lambda} \Vert^2 + \langle \eta_{\lambda}, G(u_{\lambda z}) - L z - L \lambda \rangle.
\end{equation}

Substituting (\ref{pf:forward:conv:3}) and (\ref{pf:forward:conv:4}) into (\ref{pf:forward:conv:2}), we derive
\begin{equation*}
\begin{aligned}
\dot S^f_{\lambda z} \le&
-\delta^f_{\tau_{\lambda}}\Vert \tau_{\lambda}\Vert^2 
- \delta^f_{\tau_z}\Vert \tau_z \Vert^2
+\big\langle z - z^*, L\lambda \big\rangle\\
&+ \big\langle \lambda - \lambda^*, G(u_{\lambda z}) - L \lambda - Lz \big\rangle.
\end{aligned}
\end{equation*}
Clearly, $\lambda(t) - \lambda^* \in \mathcal T_{\mathbb R^{\tilde m}_+}(\lambda^*)$, $L \lambda^* = \mathbf{0}$, and
$G(u_{\lambda z}^*) - L \lambda^* - Lz^* \in \mathcal N_{\mathbb R_+^{\tilde m}}(\lambda^*)$.
Thus, 
$\langle \lambda - \lambda^*, G(u_{\lambda z}^*) - L z^* \rangle \le 0$, and 
\begin{equation*}
\begin{aligned}
\dot S^f_{\lambda z} \le & 
- \delta^f_{\tau_{\lambda}}\Vert \tau_{\lambda} \Vert^2
- \delta^f_{\tau_z}\Vert \tau_z \Vert^2
- \lambda^T L \lambda\\
&+ \big\langle \lambda - \lambda^*, G(u_{\lambda z}) - G(u_{\lambda z}^*)\big\rangle.
\end{aligned}
\end{equation*}
Due to the convexity of $G$ and $\lambda \ge 0$,
$\big\langle \lambda, G(u_{\lambda z}) - G(u_{\lambda z}^*) \big\rangle \le 
\langle \lambda, \nabla G(u_{\lambda z})(u_{\lambda z} - u_{\lambda z}^*)\rangle 
= \langle y_{\lambda z}, u_{\lambda z} - u^*_{\lambda z}\rangle,$
and moreover,
$- \big\langle \lambda^*, G(u_{\lambda z}) - G(u_{\lambda z}^*) \big\rangle \le - \langle y_{\lambda z}^*, u_{\lambda z} - u^*_{\lambda z}\rangle.$
Therefore,
\begin{equation}
\begin{aligned}
\label{pf:forward:conv:5}
\dot S^f_{\lambda z} \le &
- \delta^f_{\tau_{\lambda}}\Vert \tau_{\lambda} \Vert^2
- \delta^f_{\tau_z}\Vert \tau_z \Vert^2 
- \lambda^T L \lambda \\
& + \langle y_{\lambda z}- y_{\lambda z}^*, u_{\lambda z} - u^*_{\lambda z}\rangle.
\end{aligned}
\end{equation}
In conclusion,
$\dot S^f_{\lambda z} \le \langle \tilde u_{\lambda z}, \tilde y_{\lambda z} \rangle$,
and part (ii) holds.

(iii) Construct a Lyapunov function candidate as 
\begin{equation}
\label{pf:forward:conv:6}
S^f = S^f_x + S^f_{\lambda z}.
\end{equation}

Note that
$u_x = -y_{\lambda z}$ and $u_{\lambda z} = y_x$.
Combining (\ref{pf:forward:conv:1}) with (\ref{pf:forward:conv:5}), we derive
\begin{equation}
\begin{aligned}
\label{pf:forward:conv:7}
\dot S^f  \le& - \delta^f_{\tau_x}\Vert \tau_x \Vert^2 
- \delta^f_{\tau_{\lambda}}\Vert \tau_{\lambda} \Vert^2
- \delta^f_{\tau_z}\Vert \tau_z \Vert^2\\
&-\big\langle x - x^*, F(x) - F(x^*) \big\rangle 
- \lambda^T L \lambda.
\end{aligned}
\end{equation}

The monotonicity of $F$ implies $\dot S^f \le 0$. The state trajectory of $\big(\rho(t), \tau(t)\big)$ is bounded since $S^f$ is radially unbounded.

Let $\mathcal R^f = \{(\rho, \tau)~|~ \dot S^f = 0\} \subset \{(\rho, \tau) ~|~\tau = \mathbf{0}, L \lambda = \mathbf{0}\}$, and 
$\mathcal M^f$ be the largest invariant subset of $\bar{\mathcal R}^f$. It follows from the LaSalle's invariance principle \cite[Theorem. 4.4]{khalil2002nonlinear} that $\big(\rho(t), \tau(t)\big) \rightarrow \mathcal M^f$ as $t \to \infty$.
Let $\big(\bar \rho(t), \bar \tau(t)\big)$ be a trajectory of (\ref{alg:GGP:forward}), and $\big(\bar x(t), \bar \lambda(t), \bar z(t)\big)$ be the corresponding output.
Since $\mathcal M^f$ is an invariant set, 
$\big(\bar \rho(t), \bar \tau(t)\big) \in \mathcal M^f$
if $\big(\bar \rho(0), \bar \tau(0)\big) \in \mathcal M^f$.
Suppose $\big(\bar \rho(t), \bar \tau(t) \big) \in \mathcal M^f$ for all $t \ge 0$.
Then $\bar \tau(t) = \mathbf{0}$, $\dot {\bar \tau}(t) = \mathbf{0}$ and $L \bar \lambda(t) = \mathbf{0}$.
Furthermore, $F(\bar x(t)) + \nabla G(\bar x(t))^T \bar \lambda(t) = \mathbf{0}$ because $\Theta^f$ is full-column rank, and as a result, $\dot {\bar\rho}_x(t) = \mathbf{0}$, i.e., $ {\bar\rho}_x(t) =  {\bar\rho}_x(0)$.
Due to $L \bar\lambda(t) = \mathbf{0}$, we have $\dot {\bar\rho}_z(t) = \mathbf{0}$, i.e.,
$\bar\rho_z(t) = \bar\rho_z(0)$.
By (\ref{alg:GGP:forward}),
$\dot {\bar \rho}_{\lambda}(t) = \Pi_{\mathbb R_+^{\tilde m}}\big[\bar \rho_{\lambda}, G(\bar x(0)) - L \bar z(0)\big]. $
If $\dot {\bar \rho}_{\lambda}(t) \not= \mathbf{0}$,
then ${\bar \rho}_{\lambda}(t) = \infty$, 
which contradicts the boundedness of ${\bar \rho}_{\lambda}(t)$. 
In summary, any $\big(\rho, \tau\big) \in \mathcal M^f$ is an equilibrium point of (\ref{alg:GGP:forward}).

Finally, we show that any trajectory of $\big(\rho(t), \tau(t)\big)$ converges to an equilibrium point of (\ref{alg:GGP:forward}).
Since $\big(\rho(t), \tau(t) \big)$ converges to $\mathcal M^f$,
 there exists a strictly increasing sequence
$t_k$ with $\lim_{k \to \infty} t_k = \infty$ such that
$\lim_{k \to \infty} \big(\rho(t_k), \tau(t_k)\big) = \big(\tilde \rho, \tilde \tau \big)$,
and moreover, $\big(\tilde \rho, \tilde \tau\big) \in \mathcal M^f$.
Consider a new Lyapunov function $\tilde S^f$ defined as (\ref{pf:forward:conv:6}) by replacing $(x^*, \lambda^*, z^*)$ with $(\tilde x,\tilde \lambda, \tilde z) = (\tilde \rho_x, \tilde \rho_{\lambda}, \tilde \rho_z)$. By a similar procedure as above for $S^f$, we obtain $\dot {\tilde S}^f \le 0$. 
For any  $\epsilon > 0$, there exists $t_T$ such that 
$\tilde S^f(t_T) < \epsilon$.
Due to $\dot {\tilde S}^f \le 0$,
$\tilde S^f(t) \le \tilde S^f(t_T) < \epsilon, \forall t \ge t_T.$
As a result, the trajectory of $\big(\rho(t), \tau(t)\big)$ converges to $\big(\tilde \rho, \tilde \tau\big)$, which is an equilibrium point of (\ref{alg:GGP:forward}).
Recalling Lemma \ref{lem:forward:KKT-EQ}, 
the output trajectory of $\big(x(t), \lambda(t), z(t)\big)$ approaches a point $(\tilde x, \tilde \lambda, \tilde z)$ satisfying the KKT conditions in (\ref{KKT-EQ}), where $\tilde x$ is a GNE of (\ref{formulation}).
This completes the proof.
$\hfill\square$

\begin{remark}
\label{rmk:forward:conv}
Theorem \ref{thm:forward:convergence} indicates that the passivity of $\Sigma_x$ (\ref{sigma:x}) and $\Sigma_{\lambda z}$ (\ref{sigma:lz}) is preserved after adding the PFCs. Furthermore, dynamics (\ref{alg:GGP:forward}) can deal with a broader class of games than (\ref{alg:GP}) since it does not require $F$ to be strictly monotone.
Intuitively speaking, as a result of the compensators,
neither $\Sigma^f_x$ nor $\Sigma^f_{\lambda z}$ 
can be passive lossless as shown in (\ref{pf:forward:conv:1}) and (\ref{pf:forward:conv:5})
due to the additional negative terms introduced by the compensators.
This leads to the asymptotic convergence of (\ref{alg:GGP:forward}).
\end{remark}

\begin{remark}
We should mention that it would be not difficult to generalize the PFC (\ref{PFC:x}) to \begin{equation}
\label{PFC:x:gen}
\dot \tau^i_x = \Phi_i^f \tau^i_x + \Theta_i^f v^i_x, ~
\eta^i_x  = \Psi_i^f  \tau^i_x + \Gamma_i^f v^i_x
\end{equation}
where $\Gamma_i^f \in \mathbb R^{n_i \times n_i}$. Similarly, (\ref{PFC:z}) and (\ref{PFC:l}) can be generalized.
If all the PFCs are strictly passive, then our main results, including Lemma \ref{lem:forward:KKT-EQ} and Theorem \ref{thm:forward:convergence}, also hold. 
\end{remark}

\subsection{Connection to Other Algorithms}

Here we show that (\ref{alg:GGP:forward}) recovers two existing methods as special cases.

\emph{Modified Primal-dual Method:}
Consider a constrained optimization problem as
\begin{equation}
\label{PD:formulation}
\min\nolimits_{x \in \mathbb R^n}~f(x), ~{\rm s.t.}~Ax = b
\end{equation}
where $f: \mathbb R^n \rightarrow \mathbb R$ is a continuous, differentiable and convex function, $A \in \mathbb R^{m \times n}$ and $b \in \mathbb R^m$. 
A modified 
primal-dual method with the derivative feedback was designed as
\begin{equation}
\label{MPD}
\left\{
\begin{aligned}
\dot x &= - \nabla f(x) - A^T y, \\
\dot \rho_y &= Ax - b, ~y = \rho_y + \dot \rho_y.
\end{aligned}
\right.
\end{equation}
Dynamics (\ref{MPD}) has been applied to distributed optimization \cite{wang2010control}, robust distributed linear programming \cite{richert2015robust}, etc.
Compared to the saddle-point dynamics in  \cite{arrow1958studies}, (\ref{MPD}) can be viewed as adding a static strictly passive PFC to the integrators involved in $\rho_y$, and hence, it is a special case of (\ref{alg:GGP:forward}).

\emph{Passivity-based Gradient-Play Dynamics for GNE Seeking:} Let $H^f(s)$, $\hat H^f(s)$ and $[\bar H^f(s)]^+$ be diagonal matrices as $H^f(s) = {\rm diag}\big\{H_k^f(s)\big\}$, 
$\hat H^f(s) = {\rm diag}\big\{\hat H_{\hat k}^f(s)\big\}$ and 
$[\bar H^f(s)]^+ = {\rm diag}\big\{[\bar H_{\bar k}^f(s)]^+\big\}$
with 
\begin{equation*}
\begin{aligned}
H_k^f(s) &= \sum\nolimits_{j = 1}^{\kappa_k} {\alpha_{kj}}/{(s + \beta_{k j})} + \gamma_k, \\
\hat H_{\hat k}^f(s) &= \sum\nolimits_{\hat j = 1}^{\hat \kappa_k} {\hat \alpha_{\hat k \hat j}}/{(s + \hat \beta_{\hat k \hat j})} + \hat \gamma_{\hat k}, \\
[\bar H_{\bar k}^f(s)]^+ &\!=\! \sum\nolimits_{\bar j = 1}^{\bar \kappa_k} \Big[{\bar \alpha_{\bar k \bar j}}/{(s + \bar \beta_{\bar k \bar j})}\Big]^+ + [\bar \gamma_{\bar k}]^+
\end{aligned}
\end{equation*}
where $\beta_{k j}, \alpha_{k j}> 0, \gamma_k \ge 0, \forall j \in \{1, \dots, \kappa_k\}, k \in \{1, \dots, n\}$,
$\hat \beta_{\hat k \hat j}, \hat \alpha_{\hat k \hat j}  > 0,  \hat \gamma_{\hat k} \ge 0, \forall  \hat j \in \{1, \dots, \hat \kappa_k\}, \hat k \in \{1, \dots, \tilde m\}$, and 
$\bar \beta_{\bar k \bar j}, \bar \alpha_{\bar k \bar j} > 0, \bar \gamma_{\bar k} \ge 0, \forall  \bar j \in \{1, \dots, \bar \kappa_k\}, \bar k \in \{1, \dots, \tilde m\}$.
For these choices, (\ref{alg:GGP:forward}) is just the dynamics proposed in \cite{li2024passivity}.

\section{OUTPUT FEEDBACK COMPENSATION}
\label{sec:back}

In this section, we design another passivity-based dynamics by  adding output feedback compensators (OFCs) to (\ref{alg:GP}).
We address its convergence, investigate its relation with existing methods, and compare it with (\ref{alg:GGP:forward}) from various aspects.

\subsection{Algorithm Design}

It follows from \cite[Theorem. 6.1]{khalil2002nonlinear} that the feedback connection of two passive systems is still passive.
There were fruitful results on designing OFCs in control systems \cite{khalil2002nonlinear, sepulchre2012constructive}.
By introducing OFCs, some high-order dynamics have been developed in \cite{gadjov2022exact, gao2023second, gao2020passivity} for NE seeking.
Motivated by the observations, we consider designing a novel dynamics to solve (\ref{formulation}) by adding appropriate OFCs to (\ref{alg:GP}).

For (\ref{alg:GP}), only passive integrators $(1/s) I_{n_i}$ are used for the evolution of $x^i$ as presented in Fig. \ref{fig:GP}.
Here we add an OFC to $(1/s) I_{n_i}$.
Let the OFC be a LTI system $H_i^b(s)$  with a minimal realization as
\begin{equation}
\label{OFC:x}
\dot \xi_x^i = \Phi^b_i \xi_x^i + \Theta^b_i x^i, ~
w_x^i = \Psi^b_i \xi_x^i + \Gamma^b_i x^i
\end{equation}
where $\Phi_i^b  \in \mathbb R^{p_i^b \times p_i^b}$ is invertible, $\Theta_i^b \in \mathbb R^{p_i^b \times n_i}$ is full-column rank, $\Psi_i^b \in \mathbb R^{n_i \times p_i^b}$ is full-row rank, and $\Gamma^b_i \in \mathbb R^{n_i \times n_i}$.
Let $\xi_x := {\rm col}\{\xi_x^i\}_{i \in \mathcal I}$,
$w_x := {\rm col}\{w_x^i\}_{i \in \mathcal I}$,
$\Phi^b := {\rm blkdiag}\{\Phi_i^b\}_{i \in \mathcal I}$,
$\Theta^b := {\rm blkdiag}\{\Theta_i^b\}_{i \in \mathcal I}$, 
$\Psi^b := {\rm blkdiag}\{\Psi_i^b\}_{i \in \mathcal I}$,
and $\Gamma^b := {\rm blkdiag}\{\Gamma_i^b\}_{i \in \mathcal I}$. Then a state-space representation for $x$ is
\begin{equation}
\label{back:x:statespace}
\left \{
\begin{aligned}
\dot x& = - F(x) - \nabla G(x)^T \lambda - w_x, \\
\dot \xi_x &= \Phi^b \xi_x + \Theta^b x, ~
w_x  = \Psi^b \xi_x + \Gamma^b x.
\end{aligned}
\right.
\end{equation}

Similarly, we introduce an OFC $\hat H_i^b(s)$ for the integrators $(1/s) I_m$ involved with $z^i$ in Fig. \ref{fig:GP}.
A minimal state-space realization of $\hat H_i^b(s)$ is
\begin{equation}
\label{OFC:z}
\dot \xi_z^i = \hat \Phi^b_i \xi_z^i + \hat \Theta^b_i z^i, ~
w_z^i = \hat \Psi^b_i \xi_z^i + \hat \Gamma^b_i z^i
\end{equation}
where $\hat \Phi_i^b  \in \mathbb R^{q_i^b \times q_i^b}$ is invertible,
$\hat \Theta_i^b \in \mathbb R^{q_i^b \times m}$ is full-column rank, $\hat \Psi_i^b \in \mathbb R^{m \times q_i^b}$ is full-row rank, and $\hat \Gamma^b_i \in \mathbb R^{m \times m}$.
Let $\xi_z := {\rm col}\{\xi_z^i\}_{i \in \mathcal I}$,
$w_z := {\rm col}\{w_z^i\}_{i \in \mathcal I}$,
$\hat \Phi^b := {\rm blkdiag}\{\hat \Phi_i^b\}_{i \in \mathcal I}$,
$\hat \Theta^b := {\rm blkdiag}\{\hat \Theta_i^b\}_{i \in \mathcal I}$, $\hat \Psi^b := {\rm blkdiag}\{\hat \Psi_i^b\}_{i \in \mathcal I}$,
and $\hat \Gamma^b := {\rm blkdiag}\{\hat \Gamma_i^b\}_{i \in \mathcal I}$.
Then a state-space representation for $z$ is
\begin{equation}
\label{back:z:statespace}
\left \{
\begin{aligned}
\dot z& = L \lambda - w_z, \\
\dot \xi_z &=\hat \Phi^b \xi_z + \hat\Theta^b z, ~
w_z  = \hat \Psi^b \xi_z + \hat \Gamma^b z.
\end{aligned}
\right.
\end{equation}

We also add an OFC for $[(1/s) I_{m}]^+$ in Fig. \ref{fig:GP}. Let the OFC be a LTI system $\bar H^b_i(s)$ with a minimal realization as
\begin{equation}
\label{OFC:l}
\dot \xi_{\lambda}^i = \bar \Phi^b_i \xi_{\lambda}^i + \bar \Theta^b_i \lambda^i, ~
w_{\lambda}^i = \bar \Psi^b_i \xi_{\lambda}^i + \bar \Gamma^b_i \lambda^i
\end{equation}
where $\bar \Phi_i^b  \in \mathbb R^{r_i^b \times r_i^b}$ is invertible, $\bar \Theta_i^b \in \mathbb R^{r_i^b \times m}$ is full-column rank, $\bar \Psi_i^b \in \mathbb R^{m \times r_i^b}$ is full-row rank, and $\bar \Gamma^b_i \in \mathbb R^{m \times m}$. Denote by $\xi_{\lambda} := {\rm col}\{\xi_{\lambda}^i\}_{i \in \mathcal I}$,
$w_{\lambda} := {\rm col}\{w_{\lambda}^i\}_{i \in \mathcal I}$,
$\bar \Phi^b := {\rm blkdiag}\{\bar \Phi_i^b\}_{i \in \mathcal I}$,
$\bar \Theta^b := {\rm blkdiag}\{\bar \Theta_i^b\}_{i \in \mathcal I}$, 
$\bar \Psi^b := {\rm blkdiag}\{\bar \Psi_i^b\}_{i \in \mathcal I}$,
and
$\bar \Gamma^b := {\rm blkdiag}\{\bar \Gamma_i^b\}_{i \in \mathcal I}$.
Then  the state-space representation for $\lambda$ is
\begin{equation}
\label{back:l:statespace}
\left \{
\begin{aligned}
\dot \lambda & = \Pi_{\mathbb R_+^{\tilde m}}\big[\lambda, 
G(x) - L z - L \lambda - w_{\lambda} \big], \\
\dot \xi_{\lambda} &= \bar \Phi^b \xi_{\lambda} + \bar \Theta^b \lambda, ~
w_{\lambda}  = \bar \Psi^b \xi_{\lambda} + \bar \Gamma^b {\lambda}.
\end{aligned}
\right.
\end{equation}

\begin{remark}
Since we put $w_{\lambda}$ inside the projection operator $\Pi_{\mathbb R_+^{\tilde m}}[\lambda, \cdot]$, it holds for (\ref{back:l:statespace}) that $\lambda(t) \in \mathbb R_+^{\tilde m}$ for all $t \ge 0$ if $\lambda(0) \in \mathbb R_+^{\tilde m}$. Therefore, different from the  nonlinear compensator (\ref{PFC:l}), the LTI system (\ref{OFC:l}) can be employed here.
\end{remark}

We make the following assumption on the OFCs.
\begin{assumption}
\label{ass:back}
\begin{enumerate}[(a)]
\item For all $i \in \mathcal I$,
$H^b_i(0) = \mathbf{0}$,
$\hat H^b_i(0) = \mathbf{0}$, and
$\bar H^b_i(0) = \mathbf{0}$.

\item For all $i \in \mathcal I$, the LTI systems (\ref{OFC:x}), (\ref{OFC:z}) and (\ref{OFC:l}) are output strictly passive with respect to storage functions $S^b_{\xi_x^i} = \frac 12 \xi_x^{i,T} P_{\xi_x^i} \xi_x^i$, 
$S^f_{\xi_z^i} = \frac 12 \xi_z^{i,T} P_{\xi_z^i} \xi_z^i$, and
$S^f_{\xi_{\lambda}^i} = \frac 12 \xi_{\lambda}^{i, T} P_{\xi_{\lambda}^i} \xi_{\lambda}^i$,
where $P_{\xi_x^i}$, $P_{\xi_z^i}$ and $P_{\xi_{\lambda}^i}$ are positive definite matrices. 
\end{enumerate}
\end{assumption}

Combining (\ref{back:x:statespace}), (\ref{back:z:statespace}) with (\ref{back:l:statespace}), we propose a novel gradient-play dynamics as
\begin{equation}
\label{alg:GGP:back}
\left \{
\begin{aligned}
\dot x& = - F(x) - \nabla G(x)^T \lambda - w_x, \\
\dot \xi_x &= \Phi^b \xi_x + \Theta^b x, ~
w_x  = \Psi^b \xi_x + \Gamma^b x, \\
\dot \lambda & = \Pi_{\mathbb R_+^{\tilde m}}\big[\lambda, 
G(x) - L z - L \lambda - w_{\lambda} \big], \\
\dot \xi_{\lambda} &= \bar \Phi^b \xi_{\lambda} + \bar \Theta^b \lambda, ~
w_{\lambda}  = \bar \Psi^b \xi_{\lambda} + \bar \Gamma^b {\lambda},\\
\dot z& =L \lambda - w_z, \\
\dot \xi_z &= \hat \Phi^b \xi_z + \hat \Theta^b z, ~
w_z  = \hat \Psi^b \xi_z + \hat \Gamma^b z
\end{aligned}
\right.
\end{equation}	
where $x(0) \in \mathbb R^n$, $\xi_x(0) \in \mathbb R^{p^b}$, $\lambda(0) \in \mathbb R_+^{\tilde m}$, $\xi_x(0) \in \mathbb R^{r^b}$, $z(0) \in \mathbb R^{\tilde m}$, $\xi_z(0) \in \mathbb R^{q^b}$, $p^b = \sum_{i \in \mathcal I} p^b_i$, $q^b = \sum_{i \in \mathcal I} q^b_i$
and $r^b = \sum_{i \in \mathcal I} r^b_i$.

The following lemma addresses the relationship between equilibria of (\ref{alg:GGP:back}) and GNEs of (\ref{formulation}).

\begin{lemma}
\label{lem:back:KKT-EQ}
Let Assumptions \ref{ass:convex}, \ref{ass:graph} and \ref{ass:back} (a) hold, and $F$ be monotone.
A profile $x^*$ is a GNE of (\ref{formulation}) if and only if there exist $(\xi_x^*, \lambda^*, \xi_{\lambda}^*, z^*, \xi_z^*)$ such that $(x^*, \xi_x^*, \lambda^*, \xi_{\lambda}^*, z^*, \xi_z^*)$ is an equilibrium point of (\ref{alg:GGP:back}).
\end{lemma}

\emph{Proof:} 
(Sufficiency) Let $(x^*, \xi_x^*, \lambda^*, \xi_{\lambda}^*, z^*, \xi_z^*)$ be an equilibrium point of (\ref{alg:GGP:back}), 
and $(w_x^*, w_{\lambda}^*, w_z^*)$ be the output.

Recalling $\dot \xi_x = \mathbf{0}$ yields
$\xi_x^* = -(\Phi^b)^{-1}\Theta^b x^*$. Consequently,
$w_x^*  = \Psi^b \xi_x^* + \Gamma^b x^* = \mathbf{0}$ since $H^b(0) = -\Psi^b (\Phi^b)^{-1} \Theta^b + \Gamma^b = \mathbf{0}$.
Similarly, we have $w_{\lambda}^* = \mathbf{0}$, and $w_z^* = \mathbf{0}$.

Due to $\dot x = \mathbf{0}$ and $\dot z = \mathbf{0}$, $F(x^*) + \nabla G(x^*)^T \lambda^* = \mathbf{0}$ and $L \lambda^* = \mathbf{0}$.
By $\dot \lambda = \mathbf{0}$, 
$G(x^*) - L z^* - L \lambda^* \in \mathcal N_{\mathbb R_+^{\tilde m}}(\lambda^*)$.

In summary, the tuple $(x^*, \lambda^*, z^*)$ satisfies the KKT conditions in (\ref{KKT-EQ}), and $x^*$ is a GNE of (\ref{formulation}) by Lemma \ref{lem:GP:KKT-EQ}.

(Necessity) Let $x^*$ be a GNE of (\ref{formulation}). In light of Lemma \ref{lem:GP:KKT-EQ}, there exists $(\lambda^*, z^*)$ such that $(x^*, \lambda^*, z^*)$ satisfies (\ref{KKT-EQ}).

Take $\xi_x^* = -(\Phi^b)^{-1}\Theta^b x^*$,
$\xi_{\lambda}^* = -(\bar\Phi^b)^{-1} \bar \Theta^b \lambda^*$ and $\xi_z^* = -(\hat \Phi^b)^{-1} \hat \Theta^b z^*$.
Then $w_x^*  = \Psi^b \xi_x^* + \Gamma^b x^* = \mathbf{0}$,
$w_{\lambda}^*  = \hat \Psi^b \xi_{\lambda}^* + \hat \Gamma^b \lambda^* = \mathbf{0}$, and 
$w_z^* = \bar \Psi^b \xi_z^* + \bar \Gamma^b z^* = \mathbf{0}$.
Furthermore, we obtain ${\rm col}\{\dot x, \dot \xi_x, \dot \lambda, \dot \xi_{\lambda}, \dot z, \dot \xi_z\} = \mathbf{0}$. Thus, $(x^*, \xi_x^*, \lambda^*, \xi_{\lambda}^*, z^*, \xi_z^*)$ is an equilibrium point of (\ref{alg:GGP:back}).
This completes the proof.
$\hfill\square$

\begin{remark}
Assumption \ref{ass:back} (a) is critical to ensure the KKT conditions (\ref{KKT-EQ}) hold at the equilibria of (\ref{alg:GGP:back}), but as a result, none of the OFCs (\ref{OFC:x}), (\ref{OFC:z}) and (\ref{OFC:l}) can be strictly passive by \cite[Lemma 6.1]{khalil2002nonlinear}.
Thus, we only assume that they are output strictly passive. 
\end{remark}

\begin{remark}
To seek NEs for merely monotone games, a Tikhonov scheme was proposed in \cite{facchinei2003finite} by adding a static OFC as $y_i = - \epsilon x_i$ to (\ref{alg:GP}), where $\epsilon > 0$ is a constant.
However, the scheme can only achieve perturbed NEs. As a comparison, (\ref{alg:GGP:back}) would reach an exact GNE if it converged to an equilibrium point, by Lemma \ref{lem:back:KKT-EQ}.
\end{remark}

\subsection{Convergence Analysis}

In this subsection, we show the convergence of (\ref{alg:GGP:back}).

\begin{figure}[htp]
\centering
\includegraphics[scale=0.4]{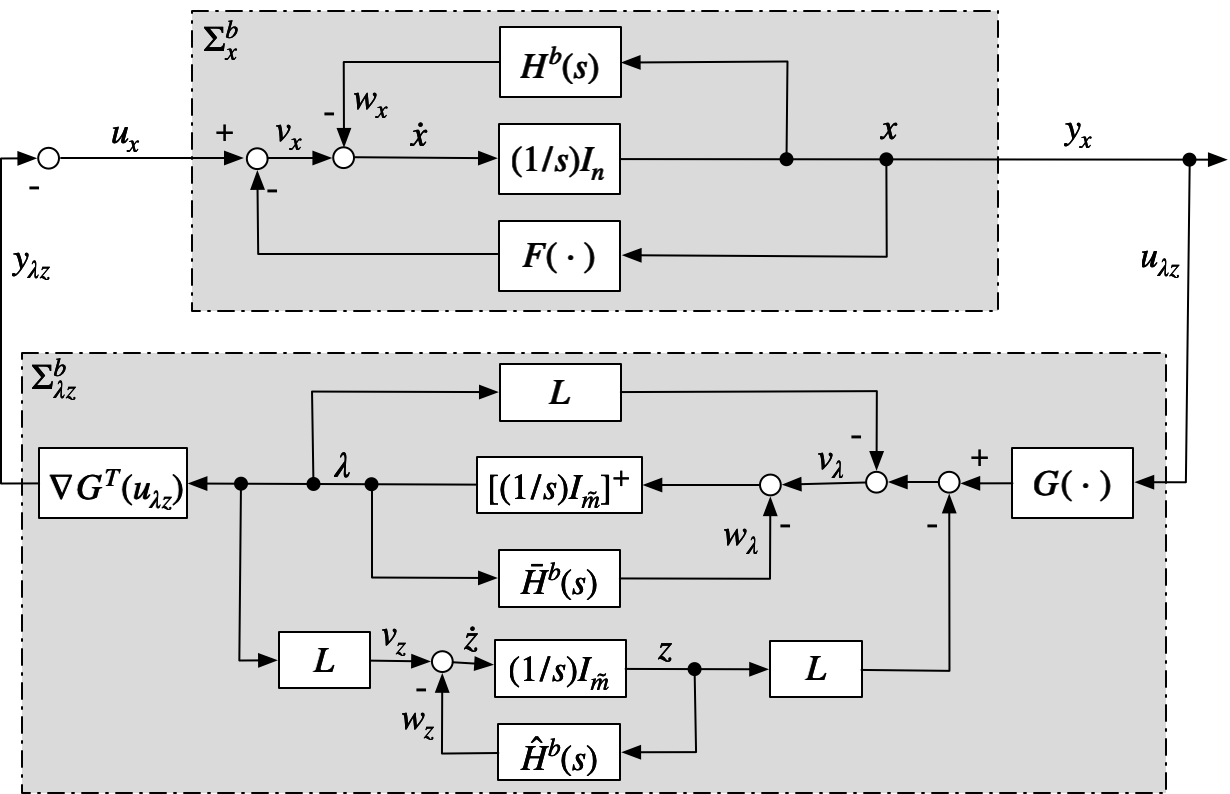}
\caption{Block diagram of dynamics (\ref{alg:GGP:back}).}
\label{fig:GGP:back}
\end{figure}

Denote by $H^b(s) := {\rm blkdiag}\{H^b_i(s)\}_{i \in \mathcal I}$,
$\bar H^b(s) := {\rm blkdiag}\{\bar H^b_i(s)\}_{i \in \mathcal I}$, and
$\hat H^b(s) := {\rm blkdiag}\{\hat H^b_i(s)\}_{i \in \mathcal I}$.
Then Fig. \ref{fig:GGP:back} shows the block diagram representation of (\ref{alg:GGP:back}). Similar to (\ref{alg:GP}),
we decompose (\ref{alg:GGP:back}) into two interconnected subsystems $\Sigma^b_x$ and $\Sigma^b_{\lambda z}$, where
\begin{equation}
\label{sigma:x:back}
\Sigma^b_x:
\left \{
\begin{aligned}
\dot x& = - F(x) + u_x - w_x, \\
\dot \xi_x &= \Phi^b \xi_x + \Theta^b x, ~
w_x  = \Psi^b \xi_x + \Gamma^b x, \\
y_x &= x
\end{aligned}
\right.
\end{equation}
and moreover,
\begin{equation}
\label{sigma:lz:back}
\Sigma^b_{\lambda z}:
\left \{
\begin{aligned}
\dot \lambda & = \Pi_{\mathbb R_+^{\tilde m}}\big[\lambda, 
G(u_{\lambda z}) - L z - L \lambda - w_{\lambda} \big], \\
\dot \xi_{\lambda} &= \bar \Phi^b \xi_{\lambda} + \bar \Theta^b \lambda, ~
w_{\lambda}  = \bar \Psi^b \xi_{\lambda} + \bar \Gamma^b {\lambda}, \\
\dot z& =L \lambda - w_z, \\
\dot \xi_z &= \hat \Phi^b \xi_z + \hat \Theta^b z, ~
w_z  = \hat \Psi^b \xi_z + \hat \Gamma^b z, \\
y_{\lambda z} &= \nabla G(u_{\lambda z})^T \lambda.
\end{aligned}
\right.
\end{equation}	

Let $\big(\xi(t), x(t), \lambda(t), z(t)\big)$ be a state trajectory of (\ref{alg:GGP:back}), and $w(t)$ be the associated output trajectory, where $\xi(t) = {\rm col}\{\xi_x(t), \xi_{\lambda}(t), 
\xi_z(t)\}$ and $w(t) = {\rm col}\{w_x(t), w_{\lambda}(t), w_z(t)\}$. Suppose $( \xi^*, x^*, \lambda^*, z^*)$ is an equilibrium point of (\ref{alg:GGP:back}), and $w^* = {\rm col}\{w_x^*, w_{\lambda}^*, w_z^*\}$ is the output, where 
$\xi^* = {\rm col}\{\xi_x^*, \xi_{\lambda}^*, \\
\xi_z^*\}$.
From the proof of Lemma \ref{lem:back:KKT-EQ},
$w^* = \mathbf{0}$.
To facilitate the analysis, we define storage functions for  
$\Sigma^b_x$ and $\Sigma^b_{\lambda z}$ as
\begin{equation}
\label{pf:back:storage:x}
S^b_x = \frac 12 \Vert x - x^*\Vert^2 + \frac 12 (\xi_x - \xi_x^*)^T P_{\xi_x} (\xi_x - \xi_x^*)
\end{equation}
and
\begin{equation}
\begin{aligned}
\label{pf:back:storage:lz}
S^b_{\lambda z} &= \frac 12 \Vert \lambda - \lambda^*\Vert^2 
+ \frac 12 (\xi_{\lambda} - \xi_{\lambda}^*)^T P_{\xi_{\lambda}} (\xi_{\lambda} - \xi_{\lambda}^*)\\
&+ \frac 12 \Vert z - z^*\Vert^2 
+ \frac 12 (\xi_z - \xi_z^*)^T P_{\xi_z} (\xi_z - \xi_z^*)
\end{aligned}
\end{equation}
where $P_{\xi_x} = {\rm blkdiag}\{P_{\xi^i_x}\}_{i \in \mathcal I}$,
$P_{\xi_{\lambda}} = {\rm blkdiag} \{P_{\xi^i_{\lambda}}\}_{i \in \mathcal I}$,
$P_{\xi_z} = {\rm blkdiag}\{P_{\xi^i_z}\}_{i \in \mathcal I}$, $P_{\xi^i_x}$, $P_{\xi_{\lambda}^i}$ and $P_{\xi^i_z}$ are given in Assumption \ref{ass:back}.
We define $\tilde u_x$, $\tilde y_x$, $\tilde u_{\lambda z}$ and $\tilde y_{\lambda z}$ by (\ref{def:notation}).

Then convergence of (\ref{alg:GGP:back}) is established as follows.

\begin{theorem}
\label{thm:back:convergence}
Consider dynamics (\ref{alg:GGP:back}). Let Assumptions \ref{ass:convex}, \ref{ass:graph} and \ref{ass:back} hold, and $F$ be monotone. 
Suppose the outputs of $H^b_i(s)$, $\hat H^b_i(s)$ and $\bar H^b_i(s)$ stay at zero only if their states are constants for all $i \in \mathcal I$.
\begin{enumerate}[(i)]
\item System $\Sigma^b_x$ (\ref{sigma:x:back}) is passive from $\tilde u_x$ to $\tilde y_x$ with respect to the storage function $S^b_x$ (\ref{pf:back:storage:x}).

\item System $\Sigma^b_{\lambda z}$ (\ref{sigma:lz:back}) is passive from $\tilde u_{\lambda z}$ to $\tilde y_{\lambda z}$ with respect to the storage function $S^b_{\lambda z}$ (\ref{pf:back:storage:lz}).

\item The trajectory of $\big(\xi(t), x(t), \lambda(t), z(t)\big)$ converges to an equilibrium point of (\ref{alg:GGP:back}), where the $x(t)$ component approaches a GNE of (\ref{formulation}).
\end{enumerate}
\end{theorem}

\emph{Proof:}
(i) Since $H_i^b$ is output strictly passive with respect to $S^b_{\xi_x^i} = \frac 12 \xi_x^{i,T} P_{\xi_x^i} \xi_x^i$, there exists $\delta^b_{\xi_x^i} > 0$ such that
$\dot S^b_{\xi_x^i} \le \langle x^i, w_x^i \rangle -\delta^b_{\xi_x^i} \Vert w_x^i \Vert^2$.
Take $\delta^b_{\xi_x} = {\rm min}\{\delta^b_{\xi^i_x} \}$.
Referring to Example $1$ in \cite{hines2011equilibrium}, we obtain
\begin{equation}
\label{pf:back:conv:1}
\dot S^b_x  \le \langle x- x^*, \dot x\rangle 
+ \langle x - x^*, w_x - w_x^*\rangle - 
\delta^b_{\xi_x} \Vert w_x - w_x^*\Vert^2.
\end{equation}
Note that $ u_x^* -F(x^*) - w_x^* = \mathbf{0}$. 
Recalling (\ref{sigma:x:back}) gives
\begin{equation}
\begin{aligned}
\label{pf:back:conv:2}
\dot S^b_x  \le& - \delta^b_{\xi_x} \Vert w_x - w_x^*\Vert^2 + \langle u_x - u_x^*, y_x - y_x^*\rangle \\
&- \big\langle  x - x^*, F(x) - F(x^*)\big\rangle.
\end{aligned}
\end{equation}

The monotonicity of $F$ implies 
$\dot S^b_x \le \langle \tilde u_x, \tilde y_x\rangle$, and hence, part (i) holds.

(ii) By a similar procedure as above for (\ref{pf:back:conv:1}), there are $\delta^b_{\xi_\lambda} > 0$ and $\delta^b_{\xi_z} > 0$ such that
\begin{equation*}
\begin{aligned}
S^b_{\lambda z} \le \langle \lambda - \lambda^*, \dot \lambda \rangle + \langle \lambda - \lambda^*, w_{\lambda} - w_{\lambda}^*\rangle 
- \delta_{\xi_\lambda} \Vert w_{\lambda} - w_{\lambda}^* \Vert^2& \\
+ \langle z - z^*, \dot z \rangle + \langle z - z^*, w_z - w_z^*\rangle - \delta_{\xi_z} \Vert w_z - w_z^* \Vert^2&.
\end{aligned}
\end{equation*}
Similar to (\ref{pf:forward:conv:3}), we have
$\langle \lambda - \lambda^*, \dot \lambda\rangle 
\le \big\langle \lambda - \lambda^*, 
G(u_{\lambda z}) - L z - L \lambda - w_{\lambda}\big\rangle$.
Note that ${\rm col}\{w_{\lambda}^*, w_z^*\} = \mathbf{0}$.
By (\ref{sigma:x:back}), 
\begin{equation*}
\begin{aligned}
\dot S^b_{\lambda z} \le& - \delta_{\xi_\lambda} \Vert w_{\lambda} - w_{\lambda}^* \Vert^2 - 
\delta_{\xi_z} \Vert w_z - w_z^* \Vert^2 \\
& + \big\langle \lambda - \lambda^*, G(u_{\lambda z}) - L z - L \lambda \big\rangle - \langle z - z^*, L \lambda\rangle.
\end{aligned}
\end{equation*}

With a similar procedure for the proof of (\ref{pf:forward:conv:5}), we obtain
\begin{equation}
\begin{aligned}
\label{pf:back:conv:3}	
\dot S^b_{\lambda z} \le& - \delta_{\xi_\lambda} \Vert w_{\lambda} - w_{\lambda}^* \Vert^2 - 
\delta_{\xi_z} \Vert w_z - w_z^* \Vert^2 \\
& + \big\langle u_{\lambda z} - u_{\lambda z}^*, y_{\lambda z} - y_{\lambda z}^* \big\rangle 
- \lambda^T L \lambda.
\end{aligned}
\end{equation}
Therefore, $\dot S^b_{\lambda z} \le \langle \tilde u_{\lambda z}, \tilde y_{\lambda z} \rangle$, and part (ii) holds.

(iii) Construct a Lyapunov function candidate as 
\begin{equation}
\label{pf:back:conv:4}
S^b = S^b_x + S^b_{\lambda z}.
\end{equation}
Combining (\ref{pf:back:conv:2}) with (\ref{pf:back:conv:3}), we have
\begin{equation*}
\begin{aligned}
\dot S^b \le& - \delta^b_{\xi_x} \Vert w_x - w_x^*\Vert^2 
- \delta_{\xi_\lambda} \Vert w_{\lambda} - w_{\lambda}^* \Vert^2 - \delta_{\xi_z} \Vert w_z - w_z^* \Vert^2 \\
&- \big\langle  x - x^*, F(x) - F(x^*) \big\rangle - 
\lambda^T L \lambda.
\end{aligned}
\end{equation*}

As a result, $\dot S^b \le 0$. Since $S^b$ is radially unbounded,
the trajectory of $\big(\xi(t), x(t), \lambda(t), z(t)\big)$ is bounded.

Let $\mathcal R^b = \{(\xi, x, \lambda, z) ~|~ \dot S^b = 0\} \subset \{(\xi, x, \lambda, z)~|~w = w^* = \mathbf{0}, L \lambda = 0\}$,
and $\mathcal M^b$ be the largest invariant subset of $\bar{\mathcal R}^b$. By the LaSalle's invariance principle \cite[Theorem. 4.4]{khalil2002nonlinear}, $\big(\xi(t), x(t), \lambda(t), z(t)\big) \rightarrow \mathcal M^b$ as $t \to \infty$.
Let $\big(\bar \xi(t), \bar x(t), \bar \lambda(t), \bar z(t)\big)$ be a trajectory of (\ref{alg:GGP:back}), and $\bar w(t)$ be the output.
Because $\mathcal M^b$ is an invariant set, 
$\big(\bar \xi(t), \bar x(t), \bar \lambda(t), \bar z(t)\big) \in \mathcal M^b$ if $\big(\bar \xi(0), \bar x(0), \bar \lambda(0), \bar z(0)\big) \in \mathcal M^b$.
Suppose that $\big(\bar \xi(t), \bar x(t), \bar \lambda(t), \bar z(t)\big) \in \mathcal M^b$ for all $t \ge 0$.
Then $\bar w(t) = w^* = \mathbf{0}$.
Since outputs of $H^b_i(s)$, $\hat H^b_i(s)$ and $\bar H^b_i(s)$ stay at zero only if their states are constants, we have $\dot{\bar \xi}(t) = \mathbf{0}$,
i.e., $\bar \xi(t) = \bar \xi(0)$.
Furthermore, by (\ref{alg:GGP:back}), 
$\bar x(t) = -(\Theta^{b,T} \Theta^b)^{-1} \Phi^b \bar \xi_x(0)$, i.e., $\dot{\bar x}(t) = \mathbf{0}$.
Similarly, ${\rm col}\{\dot {\bar \lambda}(t), \dot {\bar z}(t)\} = \mathbf{0}$. Therefore, any $(\xi, x, \lambda, z) \in \mathcal M^b$ is an equilibrium point of (\ref{alg:GGP:back}).

Finally, we show that the trajectory of $(\xi(t), x(t), \lambda(t), z(t)\big)$ converges to an equilibrium point of (\ref{alg:GGP:back}). Since $\big(\xi(t), x(t), \lambda(t), z(t)\big) \rightarrow \mathcal M^b$, there exists a strictly increasing sequence $t_k$ with $\lim_{k \to \infty} t_k = \infty$ such that
$\lim_{k \to \infty} \big(\xi(t_k), x(t_k), \lambda(t_k), z(t_k)\big) = \big(\tilde \xi, \tilde x, \tilde \lambda, \tilde z\big)$, where 
$ \big(\tilde \xi, \tilde x, \tilde \lambda, \tilde z\big) \in \mathcal M^b$.
Consider a new Lyapunov function $\tilde S^b$ defined as (\ref{pf:back:conv:4}) by replacing $(\xi^*, x^*, \lambda^*, z^*)$ with $(\tilde \xi, \tilde x,\tilde \lambda, \tilde z)$.
By a similar procedure for discussing $S^b$,
$\dot {\tilde S}^b \le 0$. For any  $\epsilon > 0$, there exists $t_T$ such that 
$\tilde S^b(t_T) < \epsilon$.
Due to $\dot {\tilde S}^b \le 0$,
$\tilde S^b(t) \le \tilde S^b(t_T) < \epsilon, \forall t \ge t_T.$
Thus, $\big(\xi(t), x(t), \lambda(t), z(t)\big)$ converges to an equilibrium point in $\mathcal M^b$.
Recalling Lemma \ref{lem:back:KKT-EQ}, the $x(t)$ component approaches a GNE of (\ref{formulation}).
This completes the proof.
$\hfill\square$

\begin{remark}
Theorem \ref{thm:back:convergence} indicates that the OFCs can enhance the convergence of (\ref{alg:GP}). It can be inferred that
$\lim_{t \to \infty} w(t) = \mathbf{0}$ by Assumption \ref{ass:back} (b).
To ensure the convergence of $\xi(t)$, an additional requirement  is imposed on the OFCs in Theorem \ref{thm:back:convergence}.
In practice, the assumption is easy to hold since OFCs are individually designed, and two instances will be given in (\ref{HA}) and (\ref{ND:OFC}).
\end{remark}

\subsection{Connection to Other Algorithms}

Here we explore the relation between (\ref{alg:GGP:back}) and two existing methods, and design a new algorithm based on (\ref{alg:GGP:back}).

\emph{Regularization Method:}
To seek an NE for (\ref{formulation}) without $X$, a regularized dynamics, given in \cite{bruck1974strongly,facchinei2003finite}, reads as
\begin{equation}
\label{RM}
\dot x = - F(x) - \delta(t) x,
\end{equation}
where $\delta(t) > 0$ is a time-varying parameter.
In fact, (\ref{RM}) is a modified gradient-play scheme by adding an OFC to (\ref{alg:GP}), where the OFC is $y= \delta(t) x$.
Different from (\ref{OFC:x}), the compensator is a time-varying system. To ensure (\ref{RM}) reaching an exact NE, it is imposed that $\lim_{t \to \infty} \delta(t) = 0$.

\emph{Heavy-Anchor Dynamics:}
To seek an NE of (\ref{formulation}), a heavy-anchor (HA) dynamics, proposed in \cite{gadjov2022exact}, is given by
\begin{equation}
\label{HA}
\left\{
\begin{aligned}
\dot x& = - F(x) - w_x, \\
\dot \xi_x &= \alpha(x - \xi_x), ~
w_x  = \beta(x- \xi_x)
\end{aligned}
\right.
\end{equation}
where $\alpha, \beta > 0$ are constants.
It is clear that (\ref{HA}) is a special case of (\ref{alg:GGP:back}) by taking $\{\Phi^b, \Theta^b, \Psi^b, \Gamma^b\} = \{-\alpha I_n, \alpha I_n, - \beta I_n, \beta I_n\}$.

Note that (\ref{alg:GGP:back}) provides ideas on how to design  algorithms in merely monotone regimes.
The following instance is given for illustration.

Consider seeking an NE for (\ref{formulation}).
Let the OFC (\ref{OFC:x}) be
\begin{equation}
\begin{aligned}
\label{ND:OFC}
\Phi_i^b \!=\! 
\begin{bmatrix}
\mathbf{0}, ~~~~I_{n_i} \\
- I_{n_i}, -I_{n_i}
\end{bmatrix}\!\!,
\Theta_i^b \!=\!
\begin{bmatrix}
\mathbf{0} \\
I_{n_i}
\end{bmatrix}\!\!,
\Psi_i^b \!=\! \begin{bmatrix}
\mathbf{0}, I_{n_i}
\end{bmatrix}\!\!,~
{\rm and}~
\Gamma_i^b \!=\! \mathbf{0}.
\end{aligned}
\end{equation}
Then an explicit algorithm is designed. It is straightforward to verify that (\ref{ND:OFC}) satisfies all the requirements in Theorem \ref{thm:back:convergence}, and hence, the algorithm can reach an exact NE of (\ref{formulation}).

\subsection{Algorithm Comparison}

Up to now, two novel dynamics (\ref{alg:GGP:forward}) and (\ref{alg:GGP:back}) have been proposed by introducing compensators to (\ref{alg:GP}). In this subsection, we compare them from various aspects.

Both  (\ref{alg:GGP:forward}) and (\ref{alg:GGP:back}) reach exact GNEs in merely monotone regimes. However, different compensators are employed in the development of the two dynamics as follows.
\begin{itemize}
\item The PFCs in  (\ref{alg:GGP:forward}) are supposed to be 
strictly passive, but the OFCs in (\ref{alg:GGP:back}) are output strictly passive.

\item  For (\ref{alg:GGP:forward}), nonlinear PFCs (\ref{PFC:l}) are necessary to ensure $\lambda(t) \in \mathbb R^{\tilde m}_+$, while linear OFCs (\ref{OFC:l}) are adopted in (\ref{alg:GGP:back}).

\item Assumption \ref{ass:back} (a) is critical to guarantee that the KKT conditions (\ref{KKT-EQ}) hold at the equilibria of (\ref{alg:GGP:back}), but a similar requirement is not imposed to (\ref{alg:GGP:forward}).
\end{itemize}

In addition to the different compensators, dynamics
(\ref{alg:GGP:forward}) and (\ref{alg:GGP:back}) are suitable for different scenarios.
\begin{itemize}
\item Consider seeking NEs for (\ref{formulation}).
Dynamics (\ref{alg:GGP:forward}) has more potential to deal with hypomonotone games than (\ref{alg:GGP:back}).
Intuitively speaking, $\Sigma_x$  (\ref{sigma:x}) is passive-short if $F$ is hypomonotone, and it can be passivated by adding input strictly passive PFCs \cite[Chapter. 2]{sepulchre2012constructive}.
The selection of PFCs can refer to \cite{li2024passivity}.

\item Dynamics (\ref{alg:GGP:back}) is easier to be extended to handle heterogeneous set constraints as opposed to (\ref{alg:GGP:forward}). For illustration, we consider seeking an NE for (\ref{formulation}) without $X$, but agent $i$ chooses its action $x^i$  from a local convex set $\Omega_i \subset \mathbb R^{n_i}$.
Introducing the OFC (\ref{OFC:x}) to dynamics (33) in \cite{gadjov2018passivity}, we develop an algorithm as
\begin{equation*}
\label{alg:GGP:back:proj}
\left\{
\begin{aligned}
\dot x^i &= \Pi_{\Omega_i}\big[x^i, - \nabla_{x^i} J_i(x^i, x^{-i}) - w_x^i \big], \\
\dot \xi_x^i &= \Phi^b_i \xi_x^i + \Theta^b_i x^i, ~
w_x^i = \Psi^b_i \xi_x^i + \Gamma^b_i x^i.
\end{aligned}
\right.
\end{equation*}
Under Assumptions \ref{ass:convex} and \ref{ass:back}, 
$x(t) = {\rm col}\{x^i(t)\}_{i \in \mathcal I}$ approaches an NE if $F$ is monotone.
As a comparison, if we generalize (\ref{alg:GGP:forward}) for the game,
the  PFC (\ref{PFC:x}) should be a nonlinear one to ensure $x_i(t) \in \Omega_i$.
Similar to (\ref{PFC:l}), the state variable of the PFC should lie in $\Omega_i$, and meanwhile, converge to $\mathbf{0}$.
This would be impossible if $\mathbf{0} \notin \Omega_i$.
\end{itemize}

\section{GENERALIZED PASSIVITY-BASED GRADIENT-PLAY DYNAMICS}
\label{sec:generalized}

In the last two sections, we designed novel dynamics by introducing passive compensators to (\ref{alg:GP}).
A natural question is how to generalize (\ref{alg:GP}) under a unifying passivity-based framework.
In this section, we develop a class of dynamics with convergence guarantees based on such a framework.

\subsection{Algorithm Design}

Referring to \cite{lessard2016analysis, zhang2021unified}, many first-order methods for smooth optimization and games, including the gradient and momentum methods, can be modeled as linear dynamical systems with nonlinear feedbacks.
Notice that the combination of the integrators $(1/s) I_{n_i}$ and the compensators  $H^f_i(s)$ (or $H^b_i(s)$) in Fig. \ref{fig:GGP:forward} (Fig. \ref{fig:GGP:back}) is also a linear system.
Thus, it is natural to generalize (\ref{alg:GP}) via substituting the passive integrators $(1/s) I_{n_i}$ in Fig. \ref{fig:GP} by a LTI system $H_i(s)$.
A minimal state-space realization of $H_i(s)$ is
\begin{equation}
\label{comp:xi:statespace}
\dot \vartheta_x^i = A_i \vartheta_x^i + B_i v_x^i, ~\\
x^i = C_i \vartheta_x^i
\end{equation}
where $A_i \in \mathbb R^{p_i \times p_i}$,
$B_i \in \mathbb R^{p_i \times n_i}$ is full-column rank, 
and $C_i \in \mathbb R^{n_i \times p_i}$ is full-row rank.
Let $\vartheta_x := {\rm col}\{\vartheta_x^i\}_{i \in \mathcal I}$,
$A := {\rm blkdiag}\{A_i\}_{i \in \mathcal I}$,
$B := {\rm blkdiag}\{B_i\}_{i \in \mathcal I}$,
and
$C := {\rm blkdiag}\{C_i\}_{i \in \mathcal I}$.
Substituting $v_x^i$ into (\ref{comp:xi:statespace}), we obtain
\begin{equation}
\label{comp:x:statespace}
\dot \vartheta_x = A\vartheta_x - B\big[F(x) + \nabla G(x)^T \lambda\big], ~
x = C \vartheta_x.
\end{equation}

Similarly, we replace the integrators $(1/s) I_m$ for $z^i$ in Fig. \ref{fig:GP} by a LTI system $\hat H_i(s)$, one of whose minimal state-space realization is 
\begin{equation}
\label{comp:zi:statespace}
\dot \vartheta_z^i = \hat A_i \vartheta_z^i + \hat B_i v_z^i, ~
z^i = \hat C_i \vartheta_z^i
\end{equation}
where $\hat A_i \in \mathbb R^{q_i \times q_i}$,
$\hat B_i \in \mathbb R^{q_i \times m}$ is full-column rank, 
and $\hat C_i \in \mathbb R^{m \times q_i}$ is full-row rank.
Let $\vartheta_z := {\rm col}\{\vartheta_z^i\}_{i \in \mathcal I}$,
$\hat A := {\rm blkdiag}\{\hat A_i\}_{i \in \mathcal I}$,
$\hat B := {\rm blkdiag}\{\hat B_i\}_{i \in \mathcal I}$,
and
$\hat C := {\rm blkdiag}\{\hat C_i\}_{i \in \mathcal I}$.
Since $v_z = {\rm col}\{v_z^i\}_{i \in \mathcal I} = L \lambda$, we have
\begin{equation}
\label{comp:z:statespace}
\dot \vartheta_z = \hat A\vartheta_z + \hat B L \lambda, ~
z = \hat C \vartheta_z.
\end{equation}

To ensure $\lambda^i(t) \in \mathbb R^m_+$, we substitute $[(1/s) I_m]^+$ in Fig. \ref{fig:GP} by a nonlinear system $[\bar H_i(s)]^+$. Similar to (\ref{def:proj:tranf}), we let a state-space realization of $[\bar H_i(s)]^+$ be
\begin{equation}
\label{comp:li:statespace}
\dot \vartheta_{\lambda}^i =
\Pi_{\mathbb R_+^{r_i}} \big[\vartheta_{\lambda}^i,
\bar A_i \vartheta_{\lambda}^i + \bar B_i v_{\lambda}^i\big], ~
\lambda^i = \max\{\mathbf{0}, \bar C_i \vartheta_{\lambda}^i\}
\end{equation}
where $\bar A_i \in \mathbb R^{r_i \times r_i}$, 
$\bar B_i \in \mathbb R^{r_i \times m}$,
and $\bar C_i \in \mathbb R^{m \times r_i}$.
We restrict $\bar A_i$, $\bar B_i$ and $\bar C_i$ to be
$\bar A_i = {\rm blkdiag}\{\bar A_{ik}\}$,
$\bar B_i = {\rm blkdiag}\{\bar B_{ik}\}$ and 
$\bar C_i = \bar B_i^T$, where
$k \in \{1, \cdots, m\}$,
$\bar A_{ik} = \\
{\rm blkdiag}\{\bar A_{ik}^{'}, \mathbf{0}\} \in \mathbb R^{r_{ik} \times r_{ik}}$,
$\bar A_{ik}^{'} \in \mathbb R^{\bar r_{ik} \times \bar r_{ik}}$ is  negative definite,
$\bar B_{ik} = [\bar b_{ik}^1, \bar b_{ik}^2, \cdots, \bar b_{ik}^{r_{ik}}]^T \in \mathbb R^{r_{ik}}$,
$\bar b_{ik}^j > 0, \forall j \in \{1, \cdots, r_{ik}\}$, $r_{ik} \ge 1$, $0 \le \bar r_{ik} < r_{ik}$, and
$\sum_{k = 1}^m r_{ik} = r_i$.

Let $\vartheta_{\lambda} := {\rm col}\{\vartheta_{\lambda}^i\}_{i \in \mathcal I}$,
$\bar A := {\rm blkdiag}\{\bar A_i\}_{i \in \mathcal I}$,
$\bar B := {\rm blkdiag}\{\bar B_i\}_{i \in \mathcal I}$,
$\bar C := {\rm blkdiag}\{\bar C_i\}_{i \in \mathcal I}$,
and $r = \sum_{i \in \mathcal I} r_i$.
Notice that $v_{\lambda} = {\rm col}\{v^i_{\lambda}\}_{i \in \mathcal I} = G(x) - L z - L \lambda$. Then
\begin{equation}
\label{comp:l:statespace}	
\dot \vartheta_{\lambda} \!=\!
\Pi_{\mathbb R_+^r} \big[\vartheta_{\lambda}, \bar A \vartheta_{\lambda} + \bar B \big(G(x) - L z - L \lambda)\big], \lambda \!=\! \max\{\mathbf{0}, \bar C \vartheta_{\lambda}\}\!.
\end{equation}

The following assumption is made.
\begin{assumption}
\label{ass:comp}
\begin{enumerate}[(a)]
\item For every $i \in \mathcal I$, $H_i(s)$ and $\hat H_i(s)$ are PR, i.e., systems (\ref{comp:xi:statespace}) and (\ref{comp:zi:statespace}) are passive with respect to storage functions $S_{\vartheta^i_x} = \frac 12 \vartheta_x^{i, T} P_{\vartheta^i_x} \vartheta_x^i$ and
$S_{\vartheta^i_z} = \frac 12 \vartheta_z^{i, T} P_{\vartheta^i_z} \vartheta_z^i$,
where $P_{\vartheta^i_x} \succ \mathbf{0}$ and $P_{\vartheta^i_z} \succ \mathbf{0}$.

\item For every $i \in \mathcal I$, system (\ref{comp:li:statespace}) is passive with respect to the storage function $S_{\vartheta^i_{\lambda}} = \frac 12 \Vert \vartheta_{\lambda}^i \Vert^2$.

\item For every $i \in \mathcal I$, there exist full-column rank matrices $\Pi_i \in \mathbb R^{p_i \times n_i}$, $\hat \Pi_i \in \mathbb R^{q_i \times m}$ and $\bar \Pi_i \in \mathbb R^{r_i \times m}$ that solve the following regulator equations: 
\begin{equation}
\begin{aligned}
\label{ass:regulator:eq}
&\mathbf{0} = A_i \Pi_i,
&\mathbf{0} = C_i \Pi_i - I_{n_i}; \\
& \mathbf{0} = \hat A_i \hat \Pi_i, 
&\mathbf{0} = \hat C_i \hat \Pi_i - I_m;\\
&\mathbf{0} = \bar A_i \bar \Pi_i, 
&\mathbf{0} = \bar C_i \bar \Pi_i - I_m,
\end{aligned}
\end{equation}
where $\bar \Pi_i$ is a nonnegative matrix.
\end{enumerate}
\end{assumption}

\begin{remark}
Assumption \ref{ass:comp} (a) and (b) are standard for the stabilization of  dynamical systems, while (c) is employed to ensure the KKT conditions (\ref{KKT-EQ}) hold at the equilibria of our generalized dynamics. 
In fact, similar assumptions have been imposed in existing literatures, like \cite{ romano2020gne}.
Referring to \cite{wen2004unifying}, we restrict $\bar A_i$, $\bar B_i$ and $\bar C_i$ to guarantee (\ref{comp:li:statespace}) satisfying (b) and (c).
\end{remark}

Combining (\ref{comp:x:statespace}), (\ref{comp:z:statespace}) with (\ref{comp:l:statespace}), our generalized gradient-play dynamics is explicitly designed as
\begin{equation}
\label{alg:GGP}
\left\{
\begin{aligned}
\dot \vartheta_x =& A \vartheta_x - B\big[F(x) + \nabla G(x)^T \lambda \big], ~
x = C \vartheta_x  \\
\dot \vartheta_z = & \hat A \vartheta_z + \hat B L \lambda, ~
z = \hat C \vartheta_z \\
\dot \vartheta_{\lambda} =& \Pi_{\mathbb R^r_+} \big[\vartheta_{\lambda}, \bar A \vartheta_{\lambda} + \bar B \big(G(x) - L z - L \lambda) \big], \\
\lambda =& {\rm max}\{\mathbf{0}, \bar C \vartheta_{\lambda}\}
\end{aligned}
\right.
\end{equation}
where $\vartheta_x(0) \in \mathbb R^{p}$, $\vartheta_z(0) \in \mathbb R^q$, $\vartheta_{\lambda}(0) \in \mathbb R_+^r$, $p = \sum_{i \in \mathcal I}p_i$ and $q = \sum_{i \in \mathcal I} q_i$.

\begin{remark}
\label{rmk:generalized:GGP}
It is clear that by a suitable choice of $H_i(s)$, $\hat H_i(s)$ and $[\bar H_i(s)]^+$, (\ref{alg:GGP}) not only captures some typical methods including the standard gradient-play scheme (\ref{alg:GP}) and HA (\ref{HA}), but also covers dynamics (\ref{alg:GGP:forward}) and (\ref{alg:GGP:back}).
\end{remark}

Let $\Pi := {\rm blkdiag}\{\Pi_i\}_{i \in \mathcal I}$,
$\hat \Pi := {\rm blkdiag}\{\hat\Pi_i\}_{i \in \mathcal I}$, and
$\bar \Pi := {\rm blkdiag}\{\bar \Pi_i\}_{i \in \mathcal I}$.
The following lemma relates the equilibria of (\ref{alg:GGP}) with GNEs of (\ref{formulation}).

\begin{lemma}
\label{lem:comp:KKT-EQ}
Suppose that Assumptions \ref{ass:convex}, \ref{ass:graph} and \ref{ass:comp} hold, and $F$ is monotone. 
\begin{enumerate}[(i)]
\item Let $(\vartheta_x^*, \vartheta_{\lambda}^*, \vartheta_z^*)$ be an equilibrium point of (\ref{alg:GGP}), and  $(x^*, \lambda^*, z^*)$ be the corresponding output. 
Then the tuple  $(x^*, \lambda^*, z^*)$ satisfies the KKT conditions (\ref{KKT-EQ}), where $x^*$ is a GNE of (\ref{formulation}).

\item If $x^*$ is a GNE of (\ref{formulation}), then there is $(\lambda^*, z^*)$ such that $(x^*, \lambda^*, z^*)$ satisfies (\ref{KKT-EQ}).
Furthermore,
$\big(\vartheta_x^*, \vartheta_{\lambda}^*,  \vartheta_z^*\big) = \big(\Pi x^*, \bar \Pi \lambda^*, \hat \Pi z^*\big)$ is an equilibrium point of (\ref{alg:GGP}).
\end{enumerate}
\end{lemma}

\emph{Proof:}
(i) Let $(\vartheta_x^{i, *}, v_x^{i, *})$ be an equilibrium pair of (\ref{comp:xi:statespace}) such that
$A_i \vartheta_x^{i, *} + B_i v_x^{i, *} = \mathbf{0}$. 

Consider $n_i < p_i$. We denote by $B_i^{\perp} \in \mathbb R^{(p_i - n_i) \times p_i}$ to be a  left-annihilator of $B_i$. Then
${\rm rank}(B_i^{\perp}) = p_i - n_i$, $B_i^{\perp} B_i = \mathbf{0}$, and 
${\rm rank} 
\begin{bmatrix}
\begin{aligned}
B_i^{\perp} \\
B_i^T
\end{aligned}	
\end{bmatrix} = p_i$.
Furthermore,
\begin{equation*}
\begin{aligned}
A_i \vartheta_x^{i, *} + B_i v_x^{i, *} = \mathbf{0} 
&\Leftrightarrow
\begin{bmatrix}
B_i^{\perp} \\
B_i^T
\end{bmatrix} 
\big(A_i \vartheta_x^{i, *} + B_i v_x^{i, *} \big)= \mathbf{0} \\
&\Leftrightarrow
\left\{
\begin{aligned}
\mathbf{0} &= B_i^{\perp} A_i \vartheta_x^{i, *},\\
\mathbf{0} &=B_i^T A_i \vartheta_x^{i, *} + B_i^T B_i v_x^{i, *}.
\end{aligned}
\right.
\end{aligned}
\end{equation*}
Thus,  $\vartheta_x^{i, *} \in {\rm Ker}(B_i^{\perp} A_i)$, and $ v_x^{i, *} = -(B_i^T B_i)^{-1} B_i^T A_i \vartheta_x^{i, *}$.
Next, we show that ${\rm Ker}(B_i^{\perp} A_i) = {\rm Ker}(A_i)$, and as a result,
$\vartheta_x^{i, *}\in {\rm Ker}(A_i)$, and 
$v_x^{i, *} = \mathbf{0}$.

Recalling (\ref{ass:regulator:eq}) gives ${\rm dim}\big({\rm Ker}(A_i)\big) \ge n_i$, and 
${\rm rank}(A_i) \le p_i - n_i$.
Since (\ref{comp:xi:statespace}) is controllable,
${\rm rank}[A_i, B_i] = p_i$, and then, ${\rm rank}(A_i) \ge p_i - n_i$ due to ${\rm rank}(B_i) = n_i$.
Consequently,
${\rm rank}(A_i) = p_i - n_i$, and
${\rm dim}\big({\rm Ker}(A_i)\big) = n_i$.
On the other hand,
${\rm rank}[A_i, B_i] = p_i$, and
${\rm rank} (B_i^{\perp} A_i) = {\rm rank}\big(B_i^{\perp} [A_i, B_i]\big) = p_i - n_i$.
Note that ${\rm dim}\big({\rm Ker}(B_i^{\perp} A_i)\big) = 
{\rm dim}({\rm Ker}(A_i)) + {\rm dim}\big({\rm Ker}(B_i^{\perp}) \cap {\rm Im}(A_i)\big) = n_i$.
It follows that 
${\rm Ker}(B_i^{\perp}) \cap {\rm Im}(A_i) = \{\mathbf{0}\}$,
and ${\rm Ker}(B_i^{\perp} A_i) = {\rm Ker}(A_i)$.
Therefore,
$\vartheta_x^{i, *} \in {\rm Ker}(A_i)$, and $v_x^{i, *} = \mathbf{0}$.

If $p_i = n_i$, then ${\rm dim}\big({\rm Ker}(A_i)\big) = p_i$
by (\ref{ass:regulator:eq}). Hence, $A_i = \mathbf{0}$,
and $v_x^{i, *} = \mathbf{0}$.

To sum up, it holds that $ \vartheta_x^{i, *} \in {\rm Ker}(A_i)$, and $v_x^{i, *} = \mathbf{0}$ for each equilibrium pair $(\vartheta_x^{i, *}, v_x^{i, *})$ of (\ref{comp:xi:statespace}) .
Thus, $\dot \vartheta_x = \mathbf{0}$ implies $F(x^*) + \nabla G(x^*)^T \lambda^* = \mathbf{0}$  for (\ref{alg:GGP}).
By a similar procedure, we conclude that $\dot \vartheta_z = \mathbf{0}$ implies $L \lambda^* = \mathbf{0}$.

Let $(\vartheta_{\lambda}^{i, *}, v_{\lambda}^{i, *})$ be an equilibrium pair of (\ref{comp:li:statespace}).
Then there exists $\zeta_4^{i, *}$ such that
\begin{equation}
\label{pf:comp:lem:ieq1}
\zeta_4^{i, *} = \bar A_i \vartheta_{\lambda}^{i, *} + \bar B_i  v_{\lambda}^{i, *} \in \mathcal N_{\mathbb R_+^r}( \vartheta_{\lambda}^{i, *}).
\end{equation}
Multiplying $\bar \vartheta_{\lambda}^{i*,T}$ to both sides of (\ref{pf:comp:lem:ieq1}), we obtain
\begin{equation}
\label{pf:comp:lem:ieq2}
\vartheta_{\lambda}^{i*,T} \bar A_i \vartheta_{\lambda}^{i, *} +  \vartheta_{\lambda}^{i*,T} \bar B_i v_{\lambda}^{i,*} = 0.
\end{equation}
It is clear that $v_{\lambda}^{i, *} \le \mathbf{0}$ due to the restrictions on $\bar A_{ik}$ and $\bar B_{ik}$ in (\ref{comp:li:statespace}).
Consequently, $\vartheta_{\lambda}^{i*,T} \bar A_i  \vartheta_{\lambda}^{i,*} \le 0$, and
$\vartheta_{\lambda}^{i*,T} \bar B_i v_{\lambda}^{i,*} \le 0$.
Furthermore, 
$\vartheta_{\lambda}^{i, *} \in {\rm Ker}(\bar A_i)$,
and $\vartheta_{\lambda}^{i *,T} \bar B_i \bar v_{\lambda}^{i, *} = 0$.
Since $\lambda^{i, *} = \bar B_i^T \bar v_{\lambda}^{i, *}$,
$\lambda^{i*, T} v_{\lambda}^{i, *}  = 0$, i.e.,
$v_{\lambda}^{i, *} \in \mathcal N_{\mathbb R_+^m}(\lambda^{i, *})$.

In conclusion, (\ref{KKT-EQ}) holds at any equilibria of (\ref{alg:GGP}). By Lemma \ref{lem:GP:KKT-EQ}, $x^*$ is a GNE of (\ref{formulation}), and part (i) is proved.

(ii) Let $x^*$ be a GNE of (\ref{formulation}). Take 
$(\lambda^*, z^*)$ such that $(x^*, \lambda^*, z^*)$ satisfies  (\ref{KKT-EQ}). Define $\vartheta_x^* = \Pi x^*$,
$\vartheta_{\lambda}^* = \bar \Pi \lambda^*$ and
$\vartheta_z^* = \hat\Pi z^*$. 
It is straightforward to verify that $\big(\vartheta_x^*, \vartheta_{\lambda}^*,  \vartheta_z^*\big)$ is an equilibrium point of (\ref{alg:GGP}).
This completes the proof.
$\hfill\square$

\begin{remark}
\label{rmk:comp:EQ}
Suppose $(x^*, \lambda^*, z^*)$ satisfies (\ref{KKT-EQ}).
Lemma \ref{lem:comp:KKT-EQ} indicates that (\ref{alg:GGP}) could achieve a GNE of (\ref{formulation}) if the trajectory of $\big(\vartheta_x(t), \vartheta_{\lambda}(t),  \vartheta_z(t)\big)$ converged to 
$\big(\Pi x^*, \bar \Pi \lambda^*, \hat \Pi z^*\big)$.
\end{remark}

\subsection{Convergence Analysis}

\begin{figure}[htp]
\centering
\includegraphics[scale=0.4]{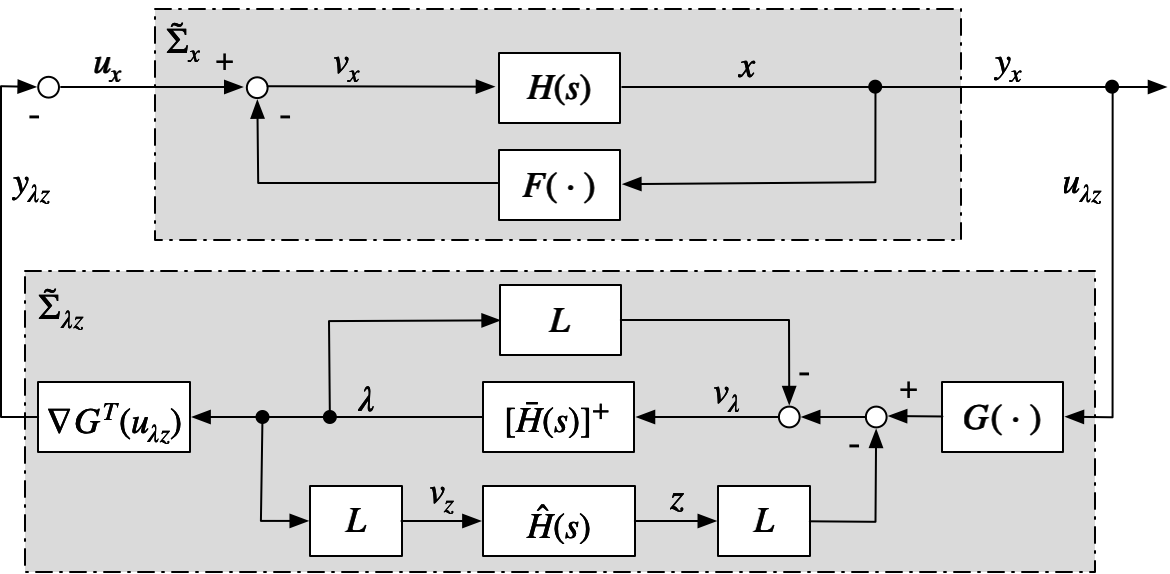}
\caption{Block diagram of dynamics (\ref{alg:GGP}).}
\label{fig:GGP:comp}
\end{figure}

Denote by
$H(s) := {\rm blkdiag}\{H_i(s)\}_{i \in \mathcal I}$,
$\hat H(s) := {\rm blkdiag}\{\hat H_i(s)\}_{i \in \mathcal I}$, and
$[\bar H(s)]^+ := {\rm blkdiag}\{[\bar H_i(s)]^+\}_{i \in \mathcal I}$.
Fig. \ref{fig:GGP:comp} shows the block diagram of (\ref{alg:GGP}). We decompose (\ref{alg:GGP}) into two interconnected systems $\tilde \Sigma_x$ and $\tilde \Sigma_{\lambda z}$, where
\begin{equation}
\label{sigma:x:comp}
\tilde \Sigma_x:
\left \{
\begin{aligned}
\dot \vartheta_x &= A \vartheta_x + B\big(u_x - F(x) \big), \\
x &= C \vartheta_x, ~
y_x = x
\end{aligned}
\right.
\end{equation}
and moreover,
\begin{equation}
\label{sigma:lz:comp}
\tilde \Sigma_{\lambda z}:
\left \{
\begin{aligned}
\dot \vartheta_{\lambda} &= \Pi_{\mathbb R^r_+} \big[ \vartheta_{\lambda}, \bar A \vartheta_{\lambda} + \bar B \big(G(u_{\lambda z}) - L z - L \lambda) \big], \\
\lambda &= {\rm max}\{\mathbf{0}, \bar C \vartheta_{\lambda}\}, \\
\dot \vartheta_z &= \hat A \vartheta_z + \hat B L \lambda, ~
z = \hat C \vartheta_z, \\
y_{\lambda z} &= \nabla G(u_{\lambda z})^T \lambda.
\end{aligned}
\right.
\end{equation}	

Suppose $x^*$ is a GNE of (\ref{formulation}). Take  $(\lambda^*, z^*)$ such that $(x^*, \lambda^*, z^*)$ satisfies (\ref{KKT-EQ}). Let $\vartheta_x^* = \Pi x^*$,
$\vartheta_{\lambda}^* = \bar \Pi \lambda^*$, and
$\vartheta_z^* = \hat \Pi z^*$.
We define storage functions for $\tilde \Sigma_x$ and $\tilde \Sigma_{\lambda z}$ as
\begin{equation}
\label{pf:comp:storage:x} 
S_x = \frac 12 (\vartheta_x - \vartheta_x^*)^T P_{\vartheta_x} (\vartheta_x - \vartheta_x^*)
\end{equation}
and 
\begin{equation}
\label{pf:comp:storage:lz}
S_{\lambda z} =   \frac 12 \Vert \vartheta_{\lambda} - \vartheta_{\lambda}^* \Vert^2
+ \frac 12 (\vartheta_z - \vartheta_z^*)^T P_{\vartheta_z}(\vartheta_z - \vartheta_z^*)
\end{equation}
where $P_{\vartheta_x} := {\rm blkdiag}\{P_{\vartheta_x^i}\}$,
$P_{\vartheta_z}  := {\rm blkdiag}\{P_{\tau^i_z}\}$,
and $P_{\vartheta_x^i}$ and $P_{\vartheta_z^i}$ are given in Assumption \ref{ass:comp}. We also define $\tilde u_x$, $\tilde y_x$, $\tilde u_{\lambda z}$ and $\tilde y_{\lambda z}$ by (\ref{def:notation}).

The following theorem establishes the convergence of (\ref{alg:GGP}).
\begin{theorem}
\label{thm:GGP:conv:strict}
Consider dynamics (\ref{alg:GGP}). Let Assumptions \ref{ass:convex}, \ref{ass:graph} and \ref{ass:comp} hold, and $F$ be monotone.
\begin{enumerate}[(i)]
\item System $\tilde \Sigma_x$ (\ref{sigma:x:comp}) is passive from $\tilde u_x$ to $\tilde y_x$ with respect the storage function $S_x$ (\ref{pf:comp:storage:x}).

\item System $\tilde \Sigma_{\lambda z}$ (\ref{sigma:lz:comp}) is passive from $\tilde u_{\lambda z}$ to $\tilde y_{\lambda z}$ with respect the storage function $S_{\lambda z}$ (\ref{pf:comp:storage:lz}).

\item The state trajectory of $\big(\vartheta_x(t), \vartheta_{\lambda}(t), \vartheta_z(t)\big)$ is bounded. If $F$ is strictly monotone, then the output trajectory of $x(t)$ converges to the GNE of (\ref{formulation}).
\end{enumerate}
\end{theorem}

\emph{Proof:}
(i) Recalling (\ref{alg:GGP}) yields
\begin{equation*}
\begin{aligned}
S_x =& \frac 12 (\vartheta_x - \vartheta_x^*)^T (A^T P_{\vartheta_x} + P_{\vartheta_x} A)(\vartheta_x - \vartheta_x^*) \\
+& (\vartheta_x - \vartheta_x^*)^T  P_{\vartheta_x} B \big[u_x - F(x) \big].
\end{aligned}
\end{equation*}
Since $H^i(s)$ is PR,  $A^T P_{\vartheta_x} + P_{\vartheta_x} A$ is negative semidefinite, and $P_{\vartheta_x} B = C^T$.
Consequently,
\begin{equation}
\label{pf:comp:conv:1}
\dot S_x \le - \big\langle F(x) - F(x^*), x - x^*\big\rangle + \langle u_x - u_x^*, y_x - y_x^* \rangle.
\end{equation}

The monotonicity of $F$ implies $\dot S_x \le \langle \tilde u_x, \tilde y_x \rangle$, and hence, part (i) holds.

(ii) It follows from (\ref{sigma:lz:comp}) that
\begin{equation*}
\begin{aligned}
\dot S_{\lambda z} =&~\frac 12 (\vartheta_z - \vartheta_z^*)^T (\hat A^T P_{\vartheta_z} + P_{\vartheta_z} \hat A)(\vartheta_z - \vartheta_z^*) \\
&+ (\vartheta_z - \vartheta_z^*)^T  P_{\vartheta_z} \hat B L \lambda + \langle \vartheta_{\lambda} - \vartheta_{\lambda}^*, \dot \vartheta_{\lambda} \rangle.
\end{aligned}
\end{equation*}
Due to the positive realness of $\hat H^i(s)$, we obtain
$\dot S_{\lambda z} \le \langle z - z^*, L\lambda \rangle
+ \langle \vartheta_{\lambda} - \vartheta_{\lambda}^*, \dot \vartheta_{\lambda} \rangle.$
Similar to (\ref{pf:forward:conv:3}), we have
\begin{equation*}
\label{pf:comp:conv:2}
\langle \vartheta_{\lambda} \!-\! \vartheta_{\lambda}^*, \dot \vartheta_{\lambda} \rangle \le \big\langle \vartheta_{\lambda} \!-\! \vartheta_{\lambda}^*, \bar A \vartheta_{\lambda}\! +\! \bar B \big[G(u_{\lambda z}) \!-\! L z\! -\! L \lambda\big] \big\rangle.
\end{equation*}

Note that $\vartheta_{\lambda}^{*} \in {\rm Ker} (\bar A)$, and
$\vartheta_{\lambda}^T \bar A \vartheta_{\lambda} \le 0$. Hence,
$\langle \vartheta_{\lambda} - \vartheta_{\lambda}^*, \dot \vartheta_{\lambda} \rangle \le \big\langle \lambda- \lambda^*, G(u_{\lambda z}) - L z - L \lambda \big\rangle$,
and moreover,
$\dot S_{\lambda z} \le \langle z - z^*, L\lambda \rangle
+ \big\langle \lambda- \lambda^*, G(u_{\lambda z}) - L z - L \lambda \big\rangle.$
By a similar procedure as the proof of 
(\ref{pf:forward:conv:5}), we have
\begin{equation}
\label{pf:comp:conv:3}
\dot S_{\lambda z} \le \langle u_{\lambda z} - u_{\lambda z}^*, y_{\lambda z} - y_{\lambda z}^*\rangle - \lambda^T L \lambda.
\end{equation}
Therefore, $\dot S_{\lambda z}  \le
\langle \tilde u_{\lambda z}, \tilde y_{\lambda z}\rangle$, and part (ii) holds.

(iii) Construct a Lyapunov function candidate as
$S = S_{x} + S_{\lambda z}$.
Combining (\ref{pf:comp:conv:1}) with (\ref{pf:comp:conv:3}), we obtain
$\dot S \le - \big\langle F(x) - F(x^*), x - x^*\big\rangle
- \lambda^T L \lambda.$
As a result, $\dot S \le 0$. Since $S$ is radially unbounded, the trajectory of $\big(\vartheta_x(t), \vartheta_{\lambda}(t), \vartheta_z(t)\big)$
is bounded.

Since $F$ is strictly monotone, $x^*$ is unique. Let $\mathcal R = \{(\vartheta_x, \vartheta_{\lambda}, \vartheta_z) ~|~ \dot S = 0\} \subset \{(\vartheta_x, \vartheta_{\lambda}, \vartheta_z)~|~x = x^*, L \lambda = \mathbf{0}\}$,
and $\mathcal M$ be the largest invariant subset of $\bar{\mathcal R}$.
Invoking the LaSalle's invariance principle \cite[Theorem. 4.4]{khalil2002nonlinear},  the trajectory of $x(t)$ converges to $x^*$, where $x^*$ is the GNE of (\ref{formulation}).
This completes the proof.
$\hfill\square$

\begin{remark}
Theorem \ref{thm:GGP:conv:strict} indicates  that the generalized gradient-play dynamics (\ref{alg:GGP}) has similar 
properties as (\ref{alg:GP}).
Since (\ref{alg:GGP}) is in a general form that covers (\ref{alg:GP}), (\ref{alg:GGP:forward}) and (\ref{alg:GGP:back}) as special cases, 
we require $F$ being strictly monotone to ensure its exact convergence.
\end{remark}

\begin{remark}
In fact, (\ref{alg:GGP}) can also be viewed from the perspective of dynamical agents. Consider $g_i(x^i) \equiv \mathbf{0}$ for explanation.
Suppose agent $i$ is associated with a LTI system 
\begin{equation}
\label{dyn:agents}
\dot \vartheta_x^i = A_i \vartheta_x^i + B_i v_x^i, 
~x^i = C_i \vartheta_x^i, 
\end{equation}
where $\vartheta_x^i$ is the state of agent $i$, $v_x^i$ is the input, and $x^i$ is the output.
We attempt to design appropriate inputs such that outputs of the agents reach an NE.
Theorem \ref{thm:GGP:conv:strict} indicates that 
$v_x^i$ can be given by $v_x^i = - \nabla_{x^i} J_i(x^i, x^{-i})$ if the system (\ref{dyn:agents}) satisfies Assumption \ref{ass:comp}. The following example is provided for illustration.
\end{remark}

\begin{example}
Consider a game played by second-order integrator agents expressed as
\begin{equation}
\label{second:order:agents}
\dot \vartheta_x^i = \begin{bmatrix}
\mathbf{0} &I_{n_i} \\
\mathbf{0} &\mathbf{0}
\end{bmatrix}\vartheta_x^i 
+  \begin{bmatrix}
\mathbf{0} \\
I_{n_i}
\end{bmatrix}
v_x^i, 
~x^i = \begin{bmatrix}
I_{n_i} ~ \mathbf{0}
\end{bmatrix} \vartheta_x^i
\end{equation}
where $\vartheta_x^i = {\rm col}\{\vartheta_{x, 1}^i, \vartheta_{x, 2}^i\} \in \mathbb R^{2n_i}$.
To render (\ref{second:order:agents}) being passive, we take
$\tilde v_i^x =  v_i^x + \frac{1}{b_i}\vartheta_{x, 2}^i$ and $\tilde x^i = \vartheta_{x, 1}^i + b_i\vartheta_{x, 2}^i$, where $b_i > 0$. Then (\ref{second:order:agents}) is cast into
\begin{equation}
\label{second:order:cast}
\dot \vartheta_x^i = \begin{bmatrix}
\mathbf{0}&I_{n_i} \\
\mathbf{0} & -\frac {1}{b_i} I_{n_i}
\end{bmatrix}\vartheta_x^i 
+  \begin{bmatrix}
\mathbf{0} \\
I_{n_i}
\end{bmatrix}
\tilde v_x^i, 
~\tilde x^i = \begin{bmatrix}
I_{n_i} ~b_i I_{n_i} 
\end{bmatrix} \vartheta_x^i.
\end{equation}
By taking $\tilde v_x^i = - \nabla_{\tilde x^i} J_i(\tilde x^i, \tilde x^{-i})$, the agents achieve an NE.
Resorting to a similar idea, the method was generalized to high-order integrator systems in \cite{romano2019dynamic}.
\end{example}


\section{PARTIAL DECISION INFORMATION SETTING}
\label{sec:partial-info}

In this section, we consider seeking a GNE of (\ref{formulation}) in a partial-decision information setting, where each agent does not have the full information on the others' actions, but only receives data from its neighbors via local communication over $\mathcal G_c$.
Similar to \cite{gadjov2018passivity, tatarenko2020geometric, bianchi2021continuous}, we let agent $i$ be endowed with an auxiliary variable $\mathbf{x}^i = {\rm col}\{\mathbf{x}^{i, j}\}_{j \in \mathcal I} \in \mathbb R^n$ that provides an estimate of all other players' strategies, where $\mathbf{x}^{i,j}$ is agent $i$'s estimate of agent $j$'s action, and $\mathbf{x}^{i, i} = x^i$ is its actual action. 
In the enlarged space, the estimate components may be different initially, but should reach consensus in the limit, i.e., $\lim_{t\to \infty} \Vert \mathbf{x}^i(t) - \mathbf{x}^j(t)\Vert = 0$.

In this case, we should notice that agent $i$ can only compute $\nabla_{x^i} J_i(\mathbf{x}^i)$ instead of the exact partial-gradient $\nabla_{x^i} J_i(x^i, x^{-i})$.
Therefore, different from $F$ in (\ref{pseudo-gradient}), we define an extended pseudogradient mapping $\mathbf{F}$ as
\begin{equation}
\label{extend:pseudo-gradient}
\mathbf{F}(\mathbf{x}) := {\rm col}\{\nabla_{x^i} J_i(\mathbf{x}^i)\}_{i \in \mathcal I} \in \mathbb R^n.
\end{equation}
Clearly, $\mathbf{F}(1_N \otimes x) = F(x)$.
Similar to (\ref{alg:GP}), we use $\bm \lambda^i \in \mathbb R^m$ to estimate $\lambda_c^*$, and introduce $\mathbf{z}^i \in \mathbb R^m$ for agent $i$. 
Let $\mathbf{x} := {\rm col}\{\mathbf{x}^i\}_{i \in \mathcal I} \in \mathbb R^{\tilde{n}}$,
$\bm \lambda := {\rm col}\{\bm \lambda^i\}_{i \in \mathcal I}\in \mathbb R^{\tilde{m}}$, $\mathbf{z} := {\rm col}\{\mathbf{z}^i\}_{i \in \mathcal I} \in \mathbb R^{\tilde{m}}$, and $\mathbf{L} = \mathcal L \otimes I_n \in \mathbb R^{\tilde n \times \tilde n}$, where $\tilde n = Nn$. Define
\begin{equation*}
\begin{aligned}
\mathcal R_i &:= 
\begin{bmatrix}
0_{n_i \times (n_{< i})}, ~~~~I_{n_i}, ~~~~0_{n_i \times (n_{> i})}
\end{bmatrix},\\
\mathcal S_i &: = \begin{bmatrix}
&I_{n_{< i}}, &0_{(n_{< i}) \times n_i}, &0_{(n_{< i}) \times (n_{> i})} \\
&0_{(n_{> i}) \times (n_{< i})}, &0_{(n_{>i}) \times n_i}, &I_{n_{> i}}
\end{bmatrix}	
\end{aligned}
\end{equation*}
where $n_{< i} = \sum_{j < i, j \in \mathcal I} n_j$, and $n_{> i} = \sum_{j > i, j \in \mathcal I} n_j$.
Let $\mathcal R := {\rm blkdiag}\{\mathcal R_i\}_{i \in \mathcal I}$ and $\mathcal S := {\rm blkdiag}\{\mathcal S_i\}_{i \in \mathcal I}$.

We make the following standard assumptions \cite{gadjov2018passivity}.
\begin{assumption}
\label{ass:partial}
\begin{enumerate}[(a)]
\item Both the pseudo-gradient mapping $F$ and the extended pseudo-gradient mapping $\mathbf{F}$ are $\theta$-Lipschitz continuous, i.e.,
$\Vert F(x) - F(y)\Vert \le \theta \Vert x - y\Vert,
\forall x, y \in \mathbb R^{n}$,
$\Vert \mathbf{F}(\mathbf{x}) - \mathbf{F}(\mathbf{y})\Vert \le \theta \Vert \mathbf{x} - \mathbf{y} \Vert,
\forall \mathbf{x}, \mathbf{y} \in \mathbb R^{\tilde{n}}$.

\item The mapping $F$ is $\mu$-strongly monotone, i.e.,
$\Vert F(x) - F(y)\Vert \ge \mu \Vert x - y\Vert,
\forall x, y \in \mathbb R^{n}$.

\item It holds that $\lambda_{\rm min}(\mathcal L) > \theta^2 / \mu + \theta$, where $\lambda_{\rm min}(\mathcal L)$ is the second minimal eigenvalue of $\mathcal L$.
\end{enumerate}	
\end{assumption}

Referring to \cite[Algorithm. 1]{bianchi2021continuous}, a fully distributed algorithm under the partial-decision information setting is proposed as
\begin{equation}
\label{alg:GP:partial}
\left\{
\begin{aligned}
\dot{\mathbf{x}} =& -{\mathcal R}^T \mathbf{F}(\mathbf{x}) - {\mathcal R}^T \nabla G(\mathcal R \mathbf{x})^T {\bm \lambda} -  \mathbf{L} \mathbf{x}, \\
\dot {\mathbf{z}} = & L {\bm \lambda}, \\
\dot {\bm \lambda} =& \Pi_{\mathbb R^{\tilde m}_+} \big[{\bm \lambda}, G(\mathcal R \mathbf{x}) - 
L \mathbf{z} - L {\bm \lambda} \big]
\end{aligned}
\right.
\end{equation}
where $\mathbf{x}(0) \in \mathbb R^{\tilde n}$, $\mathbf{z}(0) \in \mathbb R^{\tilde m}$
and ${\bm\lambda}(0) \in \mathbb R^{\tilde m}_+$.

\begin{figure}[htp]
\centering
\includegraphics[scale=0.35]{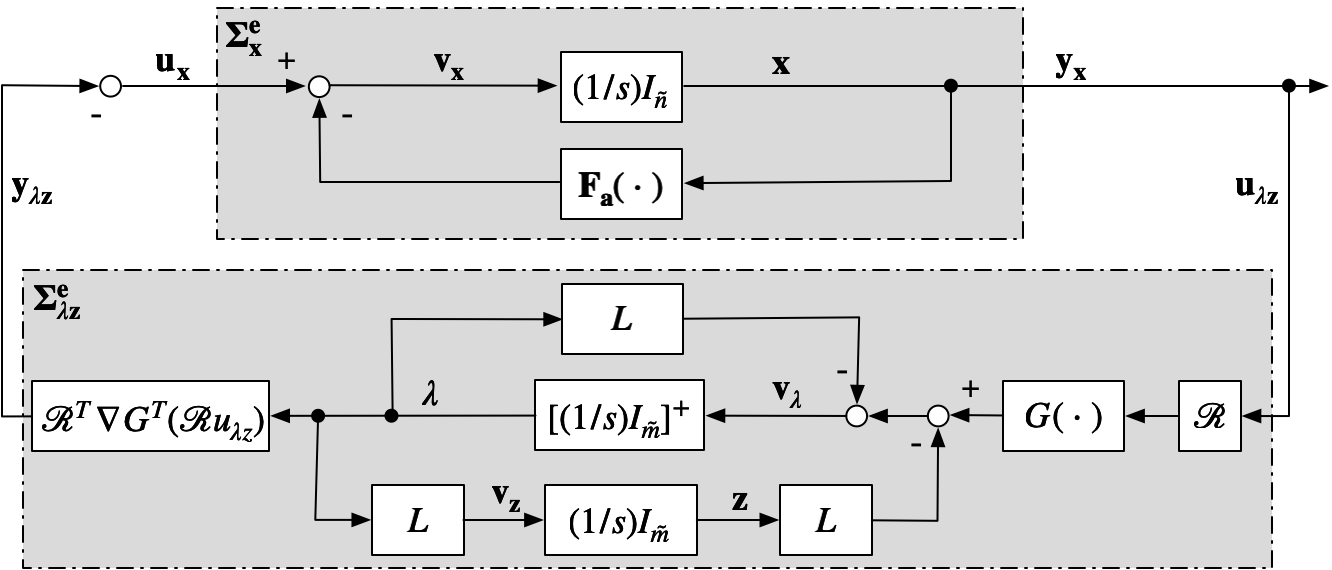}
\caption{Block diagram of dynamics (\ref{alg:GP:partial}).}
\label{fig:GP:partial}
\end{figure}

Denote by $\mathbf{F}_a(\mathbf{x}) := {\mathcal R}^T \mathbf{F}(\mathbf{x}) + \mathbf{L} \mathbf{x}$.
Then the block diagram representation of (\ref{alg:GP:partial}) is presented in Fig. \ref{fig:GP:partial}. 
We decompose (\ref{alg:GP:partial}) into two subsystems $\mathbf{\Sigma}^{\mathbf{e}}_\mathbf{x}$ and
$\mathbf{\Sigma}^{\mathbf{e}}_{{\bm \lambda} \mathbf{z}}$,
where
\begin{equation}
\label{sigma:x:partial}
\mathbf{\Sigma}^{\mathbf{e}}_\mathbf{x}: ~
\dot{\mathbf{x}} = - {\mathbf{F}}_a(\mathbf{x}) + \mathbf{u}_{\mathbf{x}}, ~
\mathbf{y}_{\mathbf{x}} = \mathbf{x}
\end{equation}
and moreover,
\begin{equation}
\label{sigma:lz:partial}
\mathbf{\Sigma}^{\mathbf{e}}_{{\bm \lambda} \mathbf{z}}: 
\left \{
\begin{aligned}
\dot {\bm\lambda} & = \Pi_{\mathbb R_+^{\tilde m}} \big[{\bm \lambda}, G(\mathcal R \mathbf{u}_{{\bm\lambda} \mathbf{z}}) - L \mathbf{z} - L {\bm \lambda} \big], \\
\dot {\mathbf{z}}& =L {\bm \lambda}, \\
\mathbf{y}_{{\bm\lambda} \mathbf{z}} &=\mathcal R^T \nabla G(\mathcal R \mathbf{u}_{{\bm\lambda} \mathbf{z}})^T {\bm \lambda}.
\end{aligned}
\right.
\end{equation}

Suppose $(\mathbf{x}^*, {\bm\lambda}^*, \mathbf{z}^*)$ is an equilibrium point of  (\ref{alg:GP:partial}).
Let $\tilde{\mathbf{u}}_{\mathbf{x}} := \mathbf{u}_{\mathbf{x}}-
\mathbf{u}_{\mathbf{x}}^*$,
$\tilde{\mathbf{u}}_{{\bm\lambda} \mathbf{z}} := \mathbf{u}_{\bm \lambda \mathbf{z}} - \mathbf{u}_{\bm \lambda \mathbf{z}}^*$,
$\tilde{\mathbf{y}}_{\mathbf{x}} := \mathbf{y}_{\mathbf{x}} - \mathbf{y}_{\mathbf{x}}^*$, and
$\tilde {\mathbf{y}}_{\bm \lambda \mathbf{z}} := \mathbf{y}_{\bm \lambda \mathbf{z}} - \mathbf{y}_{\bm \lambda \mathbf{z}}^*$,
where $\mathbf{u}_{\bm\lambda \mathbf{z}}^* \!=\! \mathbf{y}_\mathbf{x}^* \!=\! \mathbf{x}^*$ and 
$\mathbf{u}_{\mathbf{x}}^* \!=\! - \mathbf{y}_{\bm\lambda \mathbf{z}}^* \!=\! - \mathcal R^T \nabla G(\mathcal R\mathbf{x}^*)^T {\bm\lambda}^*$.
Then the convergence of (\ref{alg:GP:partial}) is addressed as follows.

\begin{theorem}
\label{thm:GP:convergence:partial}
Consider dynamics (\ref{alg:GP:partial}). Let Assumptions \ref{ass:convex}, \ref{ass:graph} and \ref{ass:partial} hold. 
\begin{enumerate}[(i)]
\item 
The mapping $\mathbf{F}_a$ is restricted strongly monotone, i.e.,
$\langle \mathbf{x} - \mathbf{x}^*, \mathbf{F}_a(\mathbf{x}) - \mathbf{F}_a(\mathbf{x}^*) \rangle \ge \delta \Vert \mathbf{x} - \mathbf{x}^*\Vert$ for some $\delta > 0$.

\item System $\mathbf{\Sigma}_{\mathbf{x}}^{\mathbf{e}}$ (\ref{sigma:x:partial}) is passive from $\mathbf{\tilde u}_{\mathbf{x}}$ to $\mathbf{\tilde y}_{\mathbf{x}}$ with respect to the storage function $\mathbf{V}_{\mathbf{x}} = \frac 12 \Vert \mathbf{x} - \mathbf{x}^*\Vert^2$.

\item System $\mathbf{\Sigma}_{{\bm \lambda} \mathbf{z}}^{\mathbf{e}}$ (\ref{sigma:lz:partial}) is passive from $\mathbf{\tilde u}_{{\bm \lambda} \mathbf{z}}$ to $\mathbf{\tilde y}_{{\bm \lambda} \mathbf{z}}$ with respect to the storage function $\mathbf{V}_{\bm \lambda \mathbf{z}} = \frac 12 \Vert \bm \lambda - \bm \lambda^*\Vert^2$ + $\frac 12 \Vert \mathbf{z} - \mathbf{z}^*\Vert^2$.

\item The trajectory $(\mathbf{x}(t), {\bm\lambda}(t), \mathbf{z}(t))$ converges to an equilibrium point of (\ref{alg:GP:partial}), where the $\mathbf{x}(t)$ component approaches $1_N \otimes x^*$, and $x^*$ is the GNE of (\ref{formulation}).
\end{enumerate}
\end{theorem}

\emph{Proof:}
Referring to Theorem $2$ in \cite{gadjov2018passivity}, part (i) holds.
By a similar procedure as that of Theorem \ref{thm:GP:convergence}, parts (\romannumeral2), (\romannumeral3) and (\romannumeral4) can be established
(see \cite{li2024passivity} for the details).
$\hfill\square$

\begin{remark}
In a partial-decision information setting, agents can only compute $\mathbf{F}$ rather than $F$ in (\ref{pseudo-gradient}).
Consequently, $F$ in Fig. \ref{fig:GP} is substituted by $\mathbf{F}_{\mathbf{a}}$ in Fig. \ref{fig:GP:partial}.
Assumption \ref{ass:partial} is made to ensure $\mathbf{F}_a$ being restricted strongly monotone, and thus,
${\bm \Sigma}_{\mathbf{x}}^\mathbf{e}$ is a passive system.
We should mention that the assumption plays a significant role on  the convergence of (\ref{alg:GP:partial}), and similar requirements have been imposed in \cite{gadjov2018passivity, bianchi2021continuous, tatarenko2020geometric}.
In the following, we generalize dynamics (\ref{alg:GP:partial}), and all the results are established under Assumption  \ref{ass:partial}.
\end{remark}

By a similar scheme presented in Sections \ref{sec:forward}, we can design a novel dynamics by adding PFCs to $(1/s) I_{\tilde n}$, $(1/s) I_{\tilde m}$ and $\big[(1/s) I_{\tilde m}\big]^+$ in Fig. \ref{fig:GP:partial}.
Let the PFCs be $\mathbf{H}^f(s) = {\rm blkdiag}\{\mathbf{H}_i^f(s)\}_{i \in \mathcal I}$,
$\hat{\mathbf{H}}^f(s) = {\rm blkdiag}\{\hat{\mathbf{H}}_i^f(s)\}_{i \in \mathcal I}$,
$[\bar{\mathbf{H}}^f(s)]^+ = {\rm blkdiag} \{[\bar{\mathbf{H}}_i^f(s)]^+\}_{i \in \mathcal I}$,
where state-space realizations of $\mathbf{H}_i^f(s)$,
$\hat{\mathbf{H}}^f(s)$ and $[\bar{\mathbf{H}}^f(s)]^+$
are similar to those of  $H^f_i(s)$, $\hat H^f_i(s)$ and $[\bar H_i^f(s)]^+$, and the only difference is that $\mathbf{H}^f_i(s)$ is a $n$-by-$n$ transfer function matrix while $H^f_i(s)$ is a $n_i$-by-$n_i$ matrix.
Suppose that $\mathbf{H}_i^f(s)$ and $\hat{\mathbf{H}}_i^f(s)$ are SPR, and the state-space realization of $[\bar{\mathbf{H}}_i^f(s)]^+$ is strictly passive.
Under Assumptions \ref{ass:convex}, \ref{ass:graph} and \ref{ass:partial}, the novel dynamics can achieve a GNE of (\ref{formulation}).
Replacing the PFCs by OFCs, we obtain another
fully distributed GNE seeking dynamics, whose convergence can be established by a similar procedure as shown in Section \ref{sec:back}.

We replace $(1/s) I_{\tilde n}$, $(1/s) I_{\tilde m}$ and $[(1/s)I_{\tilde m}]^+$ in Fig. \ref{fig:GP:partial} by passive systems $\mathbf{H}(s)$, $\hat{\mathbf{H}}(s)$
and $\big[\bar{\mathbf{H}}(s)\big]^+$. 
Let state-space realizations of $\mathbf{H}(s)$,
$\hat{\mathbf{H}}(s)$ and $[\bar{\mathbf{H}}(s)]^+$ be similar to those of $H(s)$, $\hat H(s)$ and $[\bar H(s)]$ in Section \ref{sec:generalized}.
Dynamics (\ref{alg:GP:partial}) can be generalized with convergence guarantees.
Finally, we discuss a special case of the generalization, where $g_i(x^i) \equiv \mathbf{0}$. 
Take $\mathbf{x}_{\mathbf{r}} = \mathcal R \mathbf{x}$,  
$\mathbf{x}_{\mathbf{s}} = \mathcal S \mathbf{x}$,
and $\mathbf{u}_{\mathbf{r}} = - \mathbf{F}(\mathbf{x}) - \mathcal R \mathbf{L} \mathbf{x}$.
Then $\mathbf{x} = \mathcal R^T \mathbf{x}_{\mathbf{r}}  + \mathcal S^T \mathbf{x}_{\mathbf{s}}$, and (\ref{alg:GP:partial}) reads as
$$
\dot{\mathbf{x}}_{\mathbf{r}} = \mathbf{u}_{\mathbf{r}},~
\dot{\mathbf{x}}_{\mathbf{s}} = - \mathcal S \mathbf{L} \mathbf{x}.
$$
Let $\dot{\mathbf{x}}_{\mathbf{r}} = \mathbf{u}_{\mathbf{r}}$ be substituted by
$$\dot{\bm \vartheta}_{\mathbf{r}} = \mathbf{A} \bm \vartheta_{\mathbf{r}} + \mathbf{B} \mathbf{u}_{\mathbf{r}},~
\mathbf{x}_{\mathbf{r}} = \mathbf{C} \bm \vartheta_{\mathbf{r}}$$
where $\mathbf{A} \in \mathbb R^{\tilde p \times \tilde p}$, $\mathbf{B} \in \mathbb R^{\tilde p \times n}$ and $\mathbf{C} \in \mathbb R^{n \times \tilde p}$.  Consequently,
\begin{equation}
\label{alg:GGP:partial:noconstraints}
\left \{
\begin{aligned}
\dot{\bm \vartheta}_{\mathbf{r}} &= \mathbf{A} \bm \vartheta_{\mathbf{r}} - \mathbf{B} \big( \mathbf{F}(\mathbf{x}) + \mathcal R \mathbf{L} \mathbf{x} \big), \\
\dot{\mathbf{x}}_{\mathbf{s}} &= - \mathcal S \mathbf{L} \mathbf{x}, \\
\mathbf{x} &= \mathcal R^T \mathbf{C} \bm \vartheta_{\mathbf{r}}  + \mathcal S^T \mathbf{x}_{\mathbf{s}}.
\end{aligned}
\right. 
\end{equation}
In fact, (\ref{alg:GGP:partial:noconstraints}) covers the fully distributed NE seeking methods for multi-integrator dynamical agents in \cite{romano2019dynamic}, and its convergence was analyzed in \cite{romano2020gne}.

\section{NUMERICAL SIMULATIONS}
\label{sec:simulation}

This section illustrates the applicability of our results by three examples.

\begin{example}
\label{sim:Ex1}
Consider solving Example \ref{Ex1} by (\ref{alg:GP}), (\ref{alg:GGP:forward}) and (\ref{alg:GGP:back}).
We employ two different compensators for (\ref{alg:GGP:forward}) and (\ref{alg:GGP:back}), respectively.
Fig. \ref{fig:Ex1} (a) shows the trajectories of $x(t)$, while Fig. \ref{fig:Ex1} (b) shows the trajectories of  $\log\Vert x(t) - x^*\Vert$.
The results indicate that dynamics (\ref{alg:GGP:forward}) and (\ref{alg:GGP:back}) converge to the exact NE in merely monotone regimes, whereas (\ref{alg:GP}) does not. 
\begin{figure}[htp]
\centering
\includegraphics[scale=0.22]{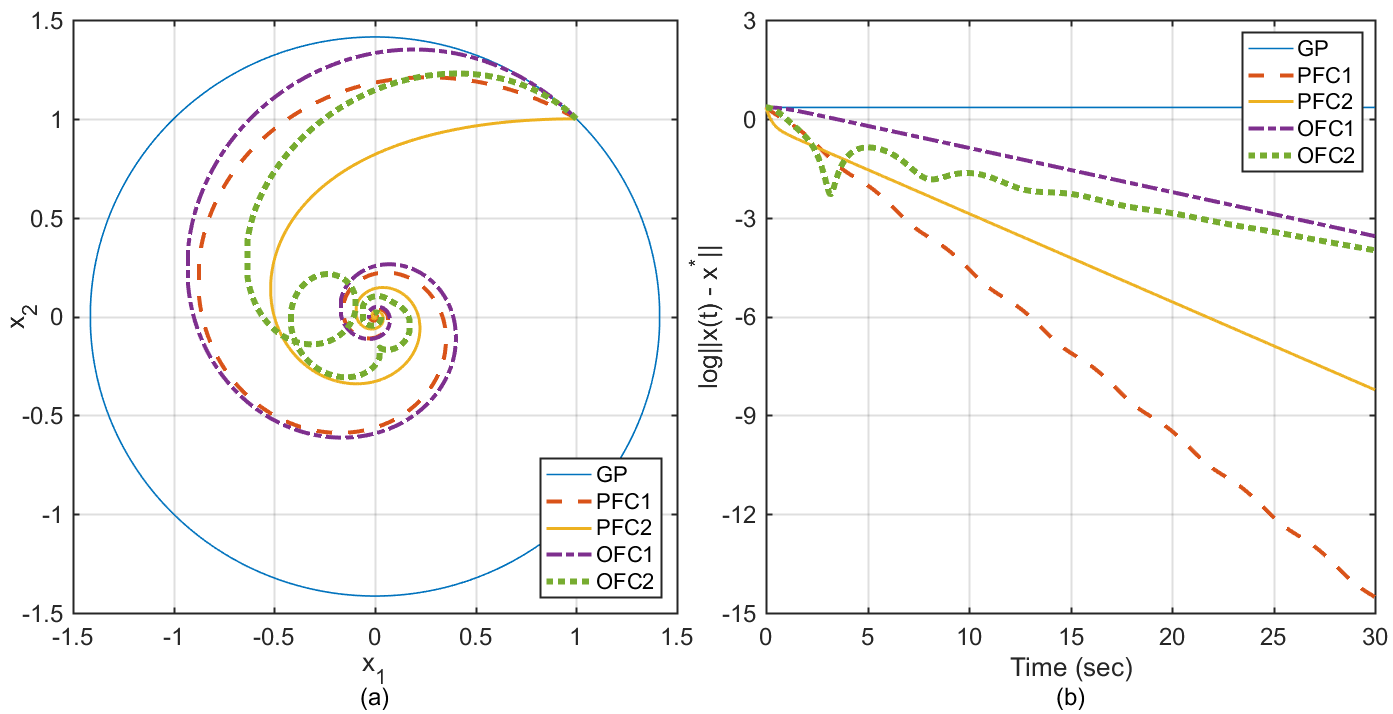}
\caption{(a) Trajectories of $x(t)$. (b) Trajectories of  $\log\Vert x(t) - x^*\Vert$. 
GP: dynamics (\ref{alg:GP});
PFC$_1$: dynamics (\ref{alg:GGP:forward}) with $H^f(s) = I_2/(s+ 1)$; 
PFC$_2$: dynamics (\ref{alg:GGP:forward}) with $H^f(s) = I_2/(s+ 4)$; 
OFC$_1$: dynamics (\ref{HA}) with $\alpha = \beta = 1$;
OFC$_2$: dynamics (\ref{alg:GGP:back}) with OFCs given by (\ref{ND:OFC}).
}
\label{fig:Ex1}
\end{figure}
\end{example}

\begin{example}
Consider a Cournot game problem as follows.
There are $N$ firms involved in the production of a homogeneous commodity that compete over $m$ markets.
Firm $i$ produces and delivers $x^i \in \mathbb R^{n_i}$ amount of products to the markets it connect with.
Let $A_i \in \mathbb R^{m \times n_i}$ (with elements $0$ or $1$) be a local matrix of firm $i$ that specifies which markets it participate in.
Define $A := [A_1, \dots, A_N] \in \mathbb R^{m \times n}$ and $x := {\rm col}\{x^i\}$. Then $Ax$ is the total product supply to all markets.
The local cost function of company $i$ is 
$$J_i(x^i, x^{-i}) = c_i(x^i) - P^T(Ax) A_i x^i$$
where $c_i(x^i) = x^{i, T} Q_i x^i + q_i^T x^i$ is the  production cost,
$P: \mathbb R^m \rightarrow \mathbb R^m$ is a price vector function given by $P(Ax) = \bar P - \Xi Ax$, $\bar P \in \mathbb R^m$ and $\Xi \in \mathbb R^{m \times m}$.
There are market capacities constraints
$\sum_{i \in \mathcal I} A_i x^i \le \sum_{i \in \mathcal I} r_i$, and the production of firm $i$ is limited by $\mathbf{0} \le x^i \le u_i$, where $r_i \in \mathbb R_+^m$ and $u_i \in \mathbb R_+^{n_i}$. 

Let $N = 5$, $m = 4$, $Q_i \succ \mathbf{0}$ with its entries randomly drawn from $(1, 4)$, each entry of $q_i$ be randomly drawn from $(0, 2)$, each entry of $\bar P$ be randomly drawn from $(10, 14)$, $\Xi \succ \mathbf{0}$ with its entries randomly drawn from $(1, 2)$, entries of $r_i$ and $u_i$ be randomly drawn from $[20, 30]$ and $[6, 14]$, respectively. 
Fig. \ref{fig:Ex2}(a) shows the trajectories of $\log(\Vert x(t) - x^*\Vert / \Vert x^*\Vert)$ under (\ref{alg:GP}), (\ref{alg:GGP:forward}) and (\ref{alg:GGP:back}), and indicates that all the three dynamics reach a GNE of the game.
In a partial-decision information setting, Fig. \ref{fig:Ex2}(b) presents the performance of  (\ref{alg:GP:partial}), and the convergence of two novel dynamics obtained by introducing PFCs and OFCs to (\ref{alg:GP:partial}).

\begin{figure}[htp]
\centering
\includegraphics[scale=0.22]{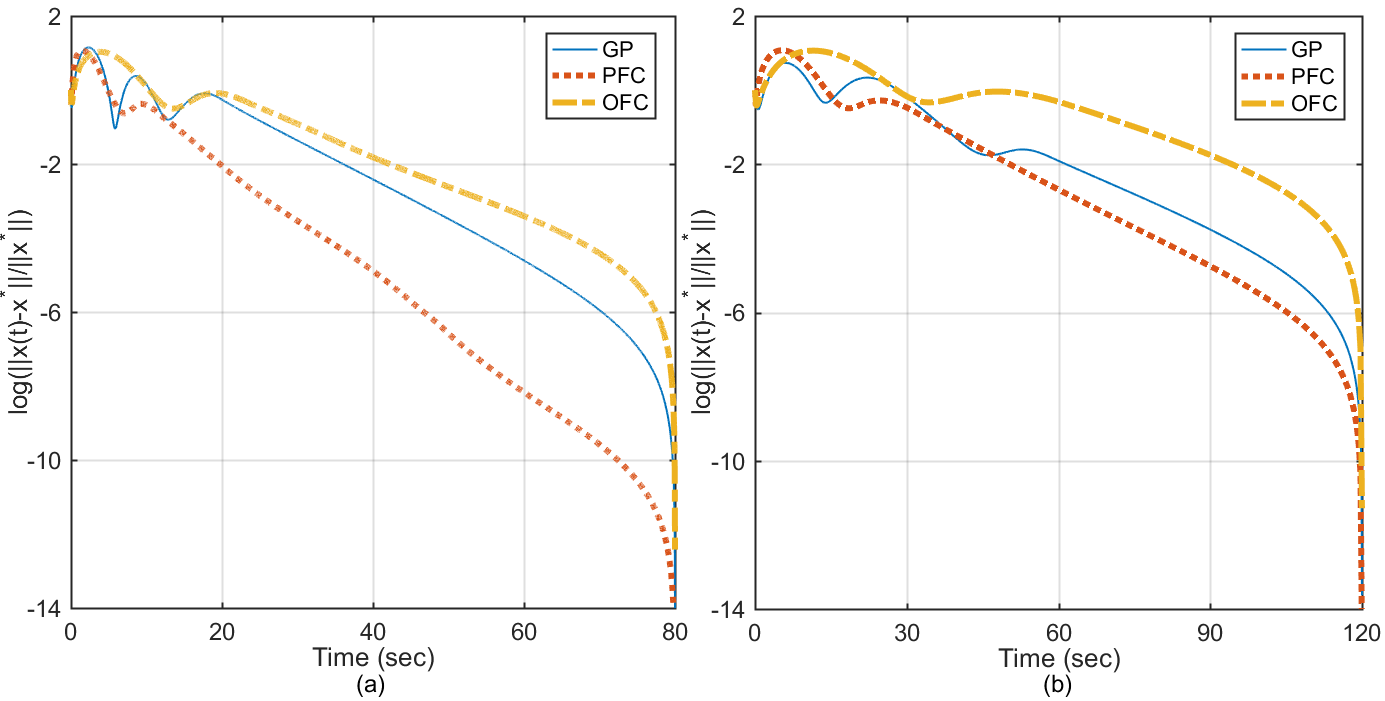}
\caption{(a) Trajectories of $\log(\Vert x(t) - x^*\Vert / \Vert x^*\Vert)$ in a full-decision information setting.
GP: dynamics (\ref{alg:GP});
PFC: dynamics (\ref{alg:GGP:forward}) with compensators given by $H_i^f(s) = I_{n_i}/(s+2)$, 
$\hat H_i^f(s) = I_{m}/(s+2)$, 
$[\bar H_i^f(s)]^+ = [I_m/(s+2)]^+$,
OFC: dynamics (\ref{alg:GGP:back}) with all compensators similar to that of HA (\ref{HA}).
(b) Trajectories of $\log(\Vert \mathbf{x}(t) - \mathbf{x}^*\Vert / \Vert \mathbf{x}^*\Vert)$ in a partial-decision information setting. GP: dynamics (\ref{alg:GP:partial});
PFC and OFC use similar compensators as those of (a).}
\label{fig:Ex2}
\end{figure}
\end{example}

\begin{example}
Consider a group of $\mathcal I = \{1, \dots, N\}$ mobile sensor devices (agents) cooperating to finish some tasks via wireless communication.
Each agent attempts to find its position in a plane to 
optimize some primary objective function, but cannot roll away from other devices.
The cost function of agent $i$ is 
\begin{equation*}
J_i(x^i, x^{-i}) = c_i(x^i) + \sum\nolimits_{j \in \mathcal I} \Vert x^i - x^j \Vert^2
\end{equation*}
where $c_i(x^i) = x^{i, T} Q_i x^i+ q_i^T x^i$, $x^i \in \mathbb R^2$, $Q_i \succ \mathbf{0}$, and $q_i \in \mathbb R^2$.
All agents send their local data to a base station, located at $\bar x$. 
To maintain acceptable levels of transmission power consumption, the average steady-state distance from agents to the base is limited by $\frac 1 N \sum_{i \in \mathcal I}\Vert x_i - \bar x \Vert^2 \le d$.

Let $N = 6$, $Q_i \succ \mathbf{0}$ with its entries randomly drawn from $(-6, 6)$, each entry of $q_i$ be randomly drawn from $(-3, 3)$, $\bar x = \mathbf{0}$, and $d = 6$.
Suppose the game is played by dynamic agents, where the dynamics of each agent is a second-order integrator (\ref{second:order:agents}). Then dynamics (\ref{alg:GGP}) is used to compute the GNE in a full-decision information setting, where a realization of $H(s)$ is (\ref{second:order:cast}), $\hat H(s) = (1/s) I_{\tilde{m}}$ and $[\bar H(s)]^+ = [(1/s) I_{\tilde{m}}]^+$.
Fig. \ref{fig:Ex3}(a) presents the trajectory of $\log(\Vert x(t) - x^* \Vert/ \Vert x^* \Vert)$, and indicates the convergence of (\ref{alg:GGP}).
Fig. \ref{fig:Ex3}(b) shows the dynamics is also convergent under a partial-decision information setting.

\begin{figure}[htp]
\centering
\includegraphics[scale=0.22]{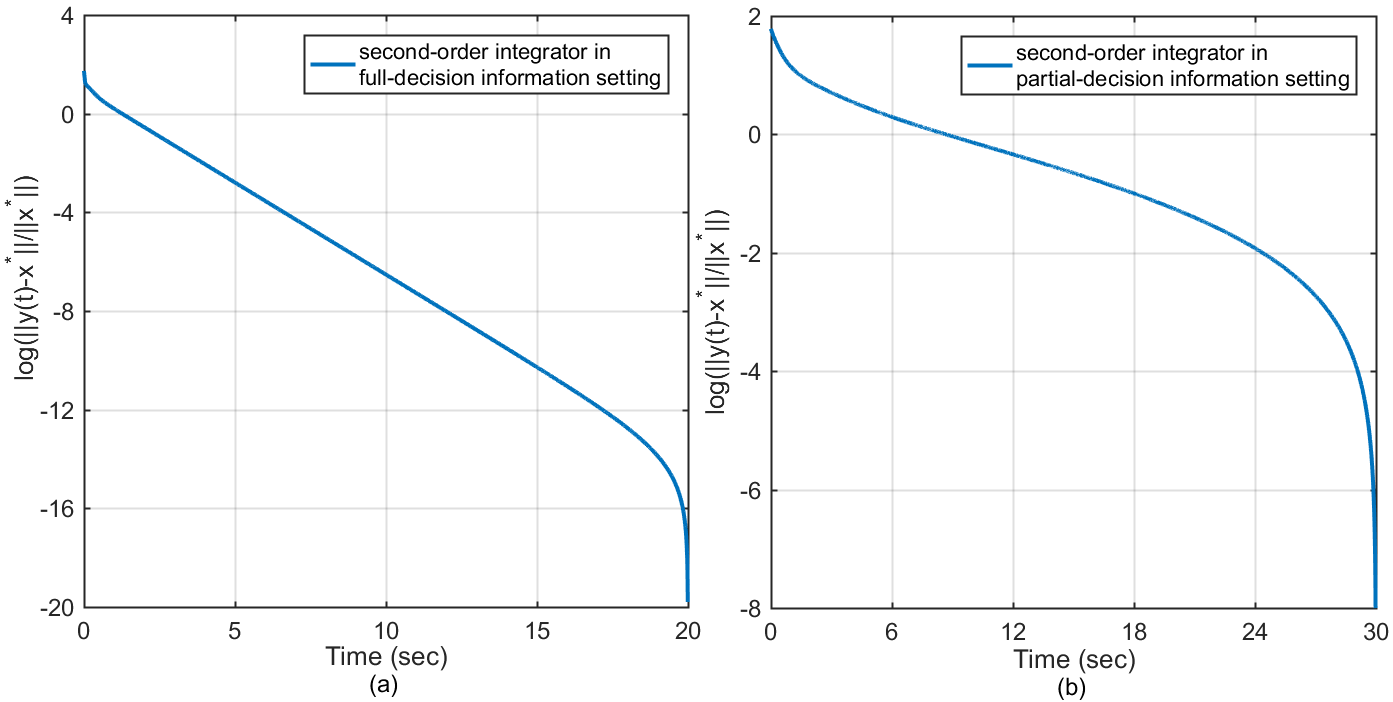}
\caption{(a) The trajectory of $\log(\Vert x(t) - x^* \Vert/ \Vert x^* \Vert)$ under  (\ref{alg:GGP}).
(b) The trajectory of $\log(\Vert \mathbf{x}(t) - \mathbf{x}^* \Vert/ \Vert \mathbf{x}^* \Vert)$ in a partial-decision information setting.}
\label{fig:Ex3}
\end{figure}
\end{example}

\section{CONCLUSION}
\label{sec:conclusion}

This paper was devoted to distributed generalized Nash equilibrium (GNE) seeking for noncooperative games with nonlinear coupled constraints. Firstly, the convergence of a typical gradient-play dynamics was analyzed based on the concept of passivity. Then two novel dynamics were proposed by parallel feedforward and output feedback compensation on the gradient-play scheme. 
Both of them allowed a relaxation of the strict monotonicity on pseudo-gradients with convergence guarantees, but they were suitable for different scenarios, i.e., one dynamics had potential to deal with hypomonotone games, while the other was powerful to handle games with heterogeneous set constraints.
After that, a class of passivity-based dynamics was further explored via substituting integrators in the gradient-play dynamics by linear-time invariant (LTI) systems.
It was established that the dynamics could reach the GNEs of strictly monotone games if the LTI systems were passive, and meanwhile, solved a group of regulator equations.
Finally, all the proposed dynamics were generalized to partial-decision information settings, and illustrative examples  were carried out for verification.

\bibliographystyle{IEEEtran}
\bibliography{references.bib}

\begin{IEEEbiography}[{\includegraphics[width=1in,height=1.25in,clip,keepaspectratio]{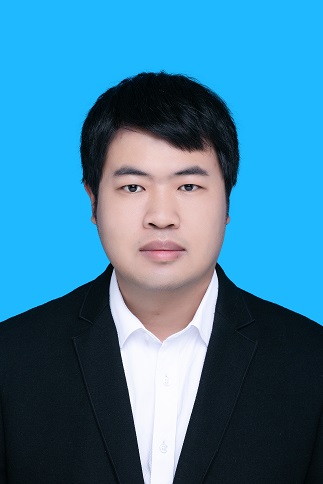}}]{Weijian Li} is a postdoctoral fellow in the Department of Electrical and Computer Engineering at University of Toronto.
He received the B.Sc. degree in mechanical engineering from Wuhan University of Technology in 2016,
and the Ph.D. degree in operations research and cybernetics from
the University of Science and Technology of China in 2021. 
He was also a joint Ph. D. student in Key Laboratory of Systems and Control, Academy of Mathematics and Systems Science, Chinese Academy of Sciences.
From Dec. 2021 to Mar. 2023, he was a researcher in the Decision Intelligence Lab, DAMO Academy, Alibaba Group.
His research interests include game theory, distributed optimization and network control.
\end{IEEEbiography}

\begin{IEEEbiography}[{\includegraphics[width=1in,height=1.25in,clip,keepaspectratio]{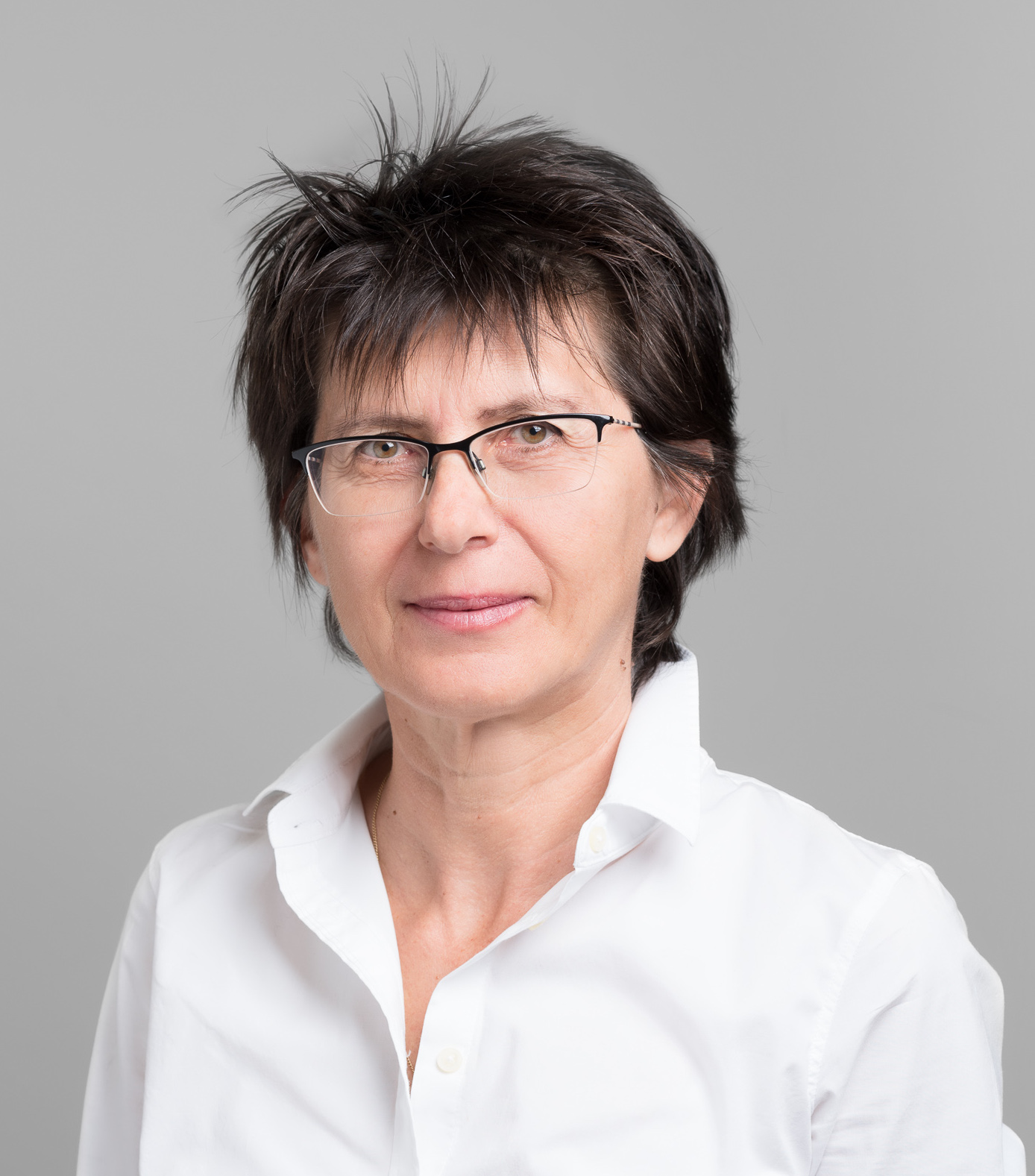}}]{Lacra Pavel} (Senior Member, IEEE) received
the Diploma of Engineering from the Technical
University of Lasi, Lasi, Romania, and the Ph.D.
degree in electrical engineering from Queen’s
University, Kingston, Canada.
After a postdoctoral stage with the National
Research Council and four years of working
in the industry, in 2002, she joined the University of Toronto, Toronto, ON, where she is
currently a Professor with the Department of
Electrical and Computer Engineering. She is the
author of the book titled Game Theory for Control of Optical Networks
(Birkhäuser/Springer, 2012). Her research interests include game theory
and distributed optimization in networks, with an emphasis on dynamics
and control aspects.
She is a Senior Editor for IEEE Transactions on Control of Networked Systems, and an Associate Editor for IEEE Transactions on Automatic Control.
\end{IEEEbiography}

\end{document}